\documentclass{amsart}
\usepackage{amsfonts}

\usepackage{amsfonts}
\usepackage{amscd}
\usepackage{amsbsy}
\usepackage{amssymb}
\usepackage{graphicx}
\usepackage[all]{xy}
\makeindex
\newtheorem{theorem}{Theorem}%[section]
\newtheorem{definition}[theorem]{Definition}

\newtheorem{lemma}[theorem]{Lemma}
\newtheorem{corollary}[theorem]{Corollary}

\newtheorem{proposition}[theorem]{Proposition}

\newtheorem{conjecture}[theorem]{Conjecture}
\newtheorem{defn-thm}[theorem]{Definition--Theorem}  %!!!!!!!!!!!!!!!!!!!!!!!!
\newtheorem{defn-prop}[theorem]{Definition--Proposition}  %!!!!!!!!!!!!!!!!!!!!!!!!

\theoremstyle{definition}
\newtheorem{oproblem}[theorem]{Open Problem}

\newtheorem{example}[theorem]{Example}

\newtheorem{remark}[theorem]{Remark}
\theoremstyle{remark}

\newcommand{\be}{\begin{equation}}
\newcommand{\ea}{\end{array}}\newcommand{\ee}[1]{\label{#1}\end{equation}}

\newcommand{\hl}{\\\hline}

\def\fract#1#2{\raise4pt\hbox{$ #1 \atop #2 $}}
\def\decdnar#1{\phantom{\hbox{$\scriptstyle{#1}$}}
\left\downarrow\vbox{\vskip15pt\hbox{$\scriptstyle{#1}$}}\right.}

\def\bfp{{\bf p}}

\def\bfz{{\bf z}}

\def\calc{{\mathcal C}}
\def\calo{{\mathcal O}}

\def\cald{{\mathcal D}}
\def\cale{{\mathcal E}}
\def\calf{{\mathcal F}}

\def\calh{{\mathcal H}}
\def\cali{{\mathcal I}}
\def\calj{{\mathcal J}}

\def\call{{\mathcal L}}
\def\calm{{\mathcal M}}

\def\calo{{\mathcal O}}

\def\calq{{\mathcal Q}}
\def\calr{{\mathcal R}}
\def\cals{{\mathcal S}}
\def\calt{{\mathcal T}}

\def\calv{{\mathcal V}}
\def\calw{{\mathcal W}}

\def\calz{{\mathcal Z}}

\def\bbc{{\mathbb C}}

\def\bbp{{\mathbb P}}

\def\bbr{{\mathbb R}}

\def\bbz{{\mathbb Z}}

\def\gra{\alpha}
\def\grb{\beta}

\def\grd{\delta}

\def\grg{\gamma}

\def\grl{\lambda}

\def\gro{\omega}

\def\grs{\sigma}
\def\grt{\tau}

\def\grG{\Gamma}

\def\grL{\Lambda}
\def\grO{\Omega}

\def\grS{\Sigma}

\def\hi{{\bar{\Phi}}}
\def\Eta{{\Upsilon}}

\def\gsp1{{\mathfrak s}{\mathfrak p}(1)}

\def\ker{\hbox{ker}}

\def\ga{{\mathfrak a}}

\def\gi{{\mathfrak i}}

\def\gn{{\mathfrak n}}

\def\gp{{\mathfrak p}}

\def\gs{{\mathfrak s}}
\def\gt{{\mathfrak t}}
\def\gu{{\mathfrak u}}

\def\gA{{\mathfrak A}}

\def\la#1{\hbox to #1pc{\leftarrowfill}}
\def\ra#1{\hbox to #1pc{\rightarrowfill}}
\def\Se{Sasaki-Einstein }

\def\BOne{{\mathchoice {\rm 1\mskip-4mu l} {\rm 1\mskip-4mu l}
                          {\rm 1\mskip-4.5mu l} {\rm 1\mskip-5mu l}}}
\def\hook{\mathbin{\hbox to 6pt{%
                 \vrule height0.4pt width5pt depth0pt
                 \kern-.4pt
                 \vrule height6pt width0.4pt depth0pt\hss}}}

\begin{document}
\bibliographystyle{amsalpha}

\title{Sasakian Geometry, Holonomy, and Supersymmetry}\footnote{During the preparation of this
work the authors
were partially supported by NSF grant DMS-0504367.}
\author{Charles P. Boyer and Krzysztof Galicki}
\address{Department of Mathematics and Statistics,
University of New Mexico, Albuquerque, NM 87131.}

\email{cboyer@math.unm.edu} \email{galicki@math.unm.edu}

%\tableofcontents

\maketitle

\tableofcontents
\bigskip
\section{Introduction}
\bigskip
Supersymmetry has emerged in physics as an attempt to unify the
way physical theories deal with bosonic and fermionic particles.
Since its birth around the early 70ties it has come to dominate
theoretical high energy physics (for a historical perspective see
\cite{KaSh00} with the introduction by Kane and Shifman, and for a
mathematical treatment see \cite{Var04}). This dominance is still
ongoing in spite of the fact that after almost 40 years there is
no single experimental evidence that would directly and
convincingly ``prove" or ``discover" the existence of
supersymmetry in nature. On the other hand, especially in the last
20 years, supersymmetry has given birth to many beautiful
mathematical theories. Gromov-Witten Theory, Seiberg-Witten
Theory, Rozansky-Witten Theory as well as the Mirror Duality
Conjecture are just a few of the more famous examples of important
and deep mathematics having its origins in the physics of various
supersymmetric theories.

Various supersymmetric field theories naturally include both
Riemannian and pseudo-Riemannian manifolds. The latter is
necessary in order to incorporate the physical space-time into the
picture while the former typically describes the geometry
associated with `invisible' extra dimensions. It is mainly in such a
context that Sasaki-Einstein manifolds appear in physics: they are
compact Einstein manifolds of positive scalar curvature that occur
in abundance in the physically interesting dimensions five and seven.
Moreover, when they are simply connected they admit real
Killing spinors. It is this last property that vitally connects
them to Supergravity, Superstring, and $M$-Theory.

The main purpose of this review is to describe geometric
properties of Sasaki-Einstein manifolds which make them
interesting in modern theoretical physics. In spite of the fact
that it is {\it supersymmetry} that connects Sasaki-Einstein
spaces to physics, it is not the purpose of this article to
describe what this concept really means to either physicists or
mathematicians. There have been many recent attempts to frame
these important notions of theoretical physics in precise
mathematical terms. This enormous task is far beyond the scope of
this article, so we refer the reader to recent monographs and
references therein \cite{QFTS99,Var04,Jos01,AJPS97}. Here we
content ourselves with providing the main theorems and results
concerning Killing spinors.

It is most remarkable that, even though Sasaki-Einstein manifolds
always have holonomy $SO(TM),$ {\it i.e.}, the holonomy of any
generic Riemannian metric, they are far from being generic. In
fact, the most interesting thing about this geometry is that it
naturally relates to several different Riemannian geometries with
reduced holonomies. It is this point that we will try to stress
throughout this article. For more detailed exposition we refer the
interested reader to our recent monograph on {\it Sasakian
Geometry} \cite{BG05}.

The key to understanding the importance of Sasakian geometry is
through its relation to K\"ahlerian geometry. Before we define
Sasakian manifolds and describe some of their elementary
properties in Section \ref{sake.snd} let us motivate things in the
more familiar context of contact and symplectic manifolds. These
two provide the mathematical foundations of Lagrangian and
Hamiltonian Mechanics. Let $(M,\eta,\xi)$ be a contact manifold
where $\eta$ is a contact form on $M$ and $\xi$ is its Reeb vector
field. It is easy to see that the cone $(C(M)=\bbr_+\times M,
\omega=d(t\eta))$ is symplectic. Likewise, the Reeb field defines
a foliation of $M$ and the transverse space $\calz$ is also
symplectic. When the foliation is regular the transverse space is
a smooth symplectic manifold giving a projection $\pi$ called
Boothby-Wang fibration, and $\pi^*\Omega=d\eta$  relates the
contact and the symplectic structures as indicated by

\def\rahmen#1#2{
\vbox{\hrule
      \hbox
        {\vrule
         \hskip#1
        \vbox{\vskip#1{}
              #2
              \vskip#1}
           \hskip#1
           \vrule}
           \hrule}}

$$\rahmen{.5cm}{
\hbox{$\begin{matrix}(\calc(M),\omega)&\hookleftarrow&(M,\eta,\xi)\\
& {}&\decdnar\pi\cr {}&{}&(\calz,\Omega).\\\end{matrix}$} }$$

We do not have any Riemannian structure yet. It is quite
reasonable to ask if there is a Riemannian metric $g$  on $M$
which ``best fits" into the above diagram. As the preferred
metrics adapted to symplectic forms are K\"ahler metrics one could
ask for the Riemannian structure which would make the cone with
the metric $\bar{g}=dt^2+t^2g$ together with the symplectic form
$\omega$ into a K\"ahler manifold. Then $\bar{g}$ and $\omega$
define a complex structure $\bar{\Phi}$. Alternatively, one could
ask for a Riemannian metric $g$ on $M$ which would define a
K\"ahler metric $h$ on $\calz$ via a Riemannian submersion.
Surprisingly, in both cases the answer to these questions leads
naturally and uniquely to {\bf Sasakian Geometry}. Our diagram
becomes

\def\rahmen#1#2{
\vbox{\hrule
      \hbox
        {\vrule
         \hskip#1
        \vbox{\vskip#1{}
              #2
              \vskip#1}
           \hskip#1
           \vrule}
           \hrule}}

           $$\rahmen{.5cm}{
\hbox{$\begin{matrix}(\calc(M),\omega,\bar{g},\bar{\Phi})&\hookleftarrow&(M,\xi,\eta,g,\Phi)\\
& {}&\decdnar\pi\cr {}&{}&(\calz,\Omega,h,J).\end{matrix}$} }
$$

From this point of view it is quite clear that K\"ahlerian and
Sasakian geometries are inseparable, Sasakian Geometry being
naturally {\it sandwiched} between two different types of
K\"ahlerian Geometry.

\medskip
\section{Cones, Holonomy, and Sasakian Geometry}\label{holonomy}
\medskip

As we have just described Sasakian manifolds can and will be
(cf. Theorem-Definition \ref{sas.def}) defined as bases of metric
cones which are K\"ahler. Let us begin with the following more
general

\begin{definition}\label{conemetric} For any Riemannian metric $g_M$ on $M,$ {\bf the warped
product
metric} on $C(M)=\bbr^+\times M$ is the Riemannian metric defined
by
$$g=dr^2+\phi^2(r)g_M\,,$$
where $r\in\bbr^+$ and $\phi=\phi(r)$ is a smooth function, called
the {\bf warping function}. If $\phi(r)=r$ then $(C(M),g)$ is
simply called the {\bf Riemannian cone} or {\bf metric cone} on
$M.$ If $\phi(r)=\sin r$ then $(C(M),g)$ is called the {\bf
sine-cone} on $M.$
\end{definition}

The relevance of sine-cones will become clear later while the
importance of metric cones in relation to the Einstein metrics can
be summarized by the following fundamental

\begin{lemma}\label{einlemma}
Let $(M,g)$ be a Riemannian manifold of dimension $n,$ and
consider $(C(M)=M\times \bbr^+,\bar{g})$ the cone on $M$ with
metric $\bar{g}=dr^2+r^2g.$ Then if $\bar{g}$ is Einstein, it is
Ricci-flat, and $\bar{g}$ is Ricci-flat if and only if $g$ is
Einstein with Einstein constant $n-1.$
\end{lemma}

Interestingly, there is a similar lemma about sine-cone metrics.

\begin{lemma}\label{sine-coneEinstein}
Let $(M^n,g)$ be an Einstein manifold with Einstein constant $n-1$
and consider $(C_s(M)=M\times (0,\pi),\bar{g}_s)$ the sine-cone on
$M$ with metric $\bar{g}_s=dr^2+(\sin^2 r)g.$ Then $\bar{g}_s$ is
Einstein with Einstein constant $n.$
\end{lemma}

It is well-known that one characterization of K\"ahlerian geometry is
via the holonomy reduction. We now recall some basic facts about the
holonomy groups of irreducible Riemannian manifolds. Let $(M,g)$ be a
Riemannian manifold and consider parallel translation defined
by the Levi-Civita connection and its associated holonomy group
which is a subgroup of the structure group $O(n,\bbr)$
($SO(n,\bbr)$ in the oriented case). Since this connection
$\nabla^g$ is uniquely associated to the metric $g$, we denote it
by ${\rm Hol}(g)$, and refer to it as the {\it Riemannian
holonomy} group or just the {\it holonomy group} when the context
is clear. Indeed, it is precisely this Riemannian holonomy that
plays an important role here. Now on a Riemannian manifold $(M,g)$
there is a canonical epimorphism $\pi_1(M)\ra{1.8} {\rm
Hol}(g)/{\rm Hol}^0(g),$ in particular, if $\pi_1(M)=0$ then ${\rm
Hol}(g)={\rm Hol}^0(g).$ In 1955 Berger proved the following
theorem \cite{Ber55} concerning Riemannian holonomy:

\begin{theorem}\label{Berger.main.thm} \index{Berger theorem} Let $(M,g)$ be an oriented
Riemannian manifold which is neither locally a Riemannian product
nor locally symmetric. Then the restricted holonomy group  ${\rm
Hol}^0(g)$ is one of the following groups listed in Table 1.1.
\end{theorem}

\begin{table}[h!]\label{bergertable}
\centering
\begin{tabular}{|l|l|l|l|}\hline
\multicolumn{4}{|c|}{{\bf Table 1.1}: Berger's Riemannian Holonomy Groups}\\
\hline
${\rm Hol}^0(g)$&dim$(M)$ & Geometry of $M$ & Comments\\
\hline \hline $SO(n)$ & $n$ &orientable Riemannian & generic
Riemannian\\ \hline $U(n)$ & $2n$ &K\"ahler & generic K\"ahler\\
\hline $SU(n)$&$2n$ & Calabi-Yau & Ricci-flat K\"ahler\\ \hline
$Sp(n)\cdot Sp(1)$ &$4n$  & quaternionic K\"ahler & Einstein\\
\hline $Sp(n)$ &$4n$ & hyperk\"ahler & Ricci-flat\\ \hline $G_2$
&$7$ & $G_2$-manifold & Ricci-flat\\ \hline
$Spin(7)$ &$8$& $Spin(7)$-manifold & Ricci-flat\\
\hline
\end{tabular}
\end{table}

Originally Berger's list included $Spin(9)$, but Alekseevsky
proved that any manifold with such holonomy must be
symmetric \cite{Al1}. In the same paper Berger also claimed a
classification of all holonomy groups of torsion-free affine
(linear) connections that act irreducibly. He produced a list of
possible holonomy representations up to what he claimed was a
finite number of exceptions. But his classification had some gaps
discovered 35 years later by Bryant \cite{Bry91}. An infinite
series of exotic holonomies was found in \cite{CMS96} and finally
the classification in the non-Riemannian affine case was completed
by Merkulov and Schwachh\"ofer \cite{MeSc99}. We refer the reader
to \cite{MeSc99} for the proof, references and the history of the
general affine case. In the Riemannian case a new geometric proof
of Berger's Theorem is now available \cite{Olm05}. An excellent
review of the subject just prior to the Merkulov and
Schwachh\"ofer's classification can be found in \cite{Bry96}. We
should add that one of the first non-trivial results concerning
manifolds with the exceptional holonomy groups of the last two
rows of Table 1.1 is due to Bonan \cite{Bon66} who established
Ricci-flatness of manifolds with parallel spinors.

Manifolds with reduced holonomy have always been very important in
physics. Partly because Calabi-Yau, hyperk\"ahler, quaternionic
K\"ahler, $G_2$ and $Spin(7)$ manifolds are automatically
Einstein. In addition, all of these spaces appear as
$\sigma$-model geometries in various supersymmetric models. What
is perhaps less known is that all of these geometries are also
related, often in more than one way, to Sasakian structures of
various flavors. Let us list all such known relations.

\begin{itemize}
\item {\bf $SO(n)$-holonomy.} As remarked this is holonomy group
of a generic metric on an oriented Riemannian manifold $(M^n,g).$
As we shall see Sasaki-Einstein metrics necessarily have maximal
holonomy. \medskip \item {\bf $U(n)$-holonomy and K\"ahler
geometry.}\index{K\"ahler manifold}
\smallskip
\begin{enumerate}
\item[(i)] Metric cone on a {\it Sasakian} manifold is {\it
K\"ahler}. \item[(ii)] Transverse geometry of a {\it Sasakian}
manifold is {\it K\"ahler}. \item[(iii)] Transverse geometry of a
{\it positive Sasakian} manifold is {\it Fano}. \item[(iv)]
Transverse geometry of a {\it Sasaki-Einstein} manifold is {\it
Fano} and {\it K\"ahler-Einstein} of {\it positive} scalar
curvature. \item[(v)] Transverse geometry of a {\it negative
Sasakian} manifold is {\it canonical}  in the sense that the
transverse canonical bundle is ample. \item[(vi)] Transverse
geometry of a {\it $3$-Sasakian} manifold is a {\it
K\"ahler-Einstein}\index{K\"ahler-Einstein manifold} with a {\it
complex contact structure}, {\it i.e.,} {\it twistor
geometry}\index{twistor geometry}.
\end{enumerate}
\medskip
\item {\bf $SU(n)$-holonomy and Calabi-Yau
geometry\index{Calabi-Yau manifold}.}
\smallskip
\begin{enumerate}
\item[(i)] Metric cone on a {\it Sasaki-Einstein} manifold is {\it
Calabi-Yau}.
\item[(ii)] Transverse geometry of a {\it null
Sasakian} manifold is {\it Calabi-Yau}.
\end{enumerate}
\medskip
\item {\bf $Sp(n)Sp(1)$-holonomy and Quaternionic K\"ahler
geometry\index{quaternionic K\"ahler manifold}.}\smallskip
\begin{enumerate}
\item[(i)] Transverse geometry of the 3-dimensional foliation of a
{\it $3$-Sasakian} manifold is {\it quaternionic-K\"ahler} of {\it
positive} scalar curvature. \item[(ii)] {\it
$3$-Sasakian}\index{$3$-Sasakian manifold} manifolds occur as
conformal infinities of complete {\it quaternionic K\"ahler}
manifolds of {\it negative} scalar curvature.
\end{enumerate}
\medskip
\item {\bf $Sp(n)$-holonomy and hyper\"ahler
geometry\index{hyperk\"ahler manifold}.}\smallskip
\begin{enumerate}
\item[(i)] Metric cone on a {\it $3$-Sasakian} manifold is {\it
hyperk\"ahler}. \item[(ii)] Transverse geometry of a {\it
null Sasakian} manifold with some additional structure is {\it
hyperk\"ahler}.
\end{enumerate}
\medskip
\item {\bf $G_2$-holonomy.}\smallskip
\begin{enumerate}
\item[(i)] The {\it `squashed' twistor space of a $3$-Sasakian}
$7$-manifold is {\it nearly K\"ahler}; hence, the metric cone on
it has holonomy inside $G_2$\index{$G_2$-holonomy manifold}.
\item[(ii)] {\it Sine-cone} on a Sasaki-Einstein $5$-manifold is
{\it nearly K\"ahler}; hence, its {\it metric cone} has holonomy
inside $G_2$.
\end{enumerate}
\medskip
\item {\bf $Spin(7)$-holonomy.}\smallskip
\begin{enumerate}
\item[(i)] The {\it `squashed' $3$-Sasakian} $7$-manifold has a 
nearly parallel $G_2$-structure; hence,
its metric cone has holonomy in $Spin(7)$. \item[(ii)] {\it
Sine-cone} on a {\it `squashed' twistor space of a $3$-Sasakian}
$7$-manifold has a {\it nearly parallel $G_2$ structure}; hence, its {\it metric
cone} has holonomy inside $Spin(7)$. \item[(iii)] {\it Sine cone
on a sine cone} on a $5$-dimensional Sasaki-Einstein base has a {\it
nearly parallel $G_2$-structure}; hence, its {\it metric cone} has
holonomy inside $Spin(7).$\index{$Spin(7)$-holonomy manifold}
\end{enumerate}
\medskip
\end{itemize}

Note that Sasakian manifolds are related to various other
geometries in two very distinct ways. On one hand we can take a
Sasakian (Sasaki-Einstein, $3$-Sasakian, etc.) manifold and
consider its metric or sine-cone. These cones frequently have
interesting geometric properties and reduced holonomy. On the
other hand, a Sasakian manifold is always naturally foliated by
one-dimensional leaves (three-dimensional leaves in addition to
the one-dimensional canonical foliation when the manifold is
$3$-Sasakian) and we can equally well consider the transverse
geometries associated to such fundamental foliations. These too
have remarkable geometric properties including reduced holonomy.
In particular, Sasakian manifolds are not just related to all of
the geometries on Berger's holonomy list, but more importantly,
they provide a {\it bridge} between the different geometries listed
there. We will investigate some of these bridges in the next two
sections.

\medskip
\section{Sasakian and K\"ahlerian geometry}\label{sake.snd}
\medskip

\begin{definition}\label{contactdef}\index{contact structure}\index{contact manifold}
A $(2n+1)$-dimensional manifold $M$ is a {\bf contact manifold} if
there exists a $1$-form $\eta$, called a {\bf contact $1$-form},
on $M$ such that
$$\eta\wedge (d\eta)^n \neq 0$$
everywhere on $M.$ A {\bf contact structure} on $M$ is an
equivalence class of such $1$-forms, where $\eta'\sim \eta$ if
there is a nowhere vanishing function $f$ on $M$ such that
$\eta'=f\eta.$
\end{definition}

\begin{lemma}\index{Reeb vector field}
On a contact manifold $(M,\eta)$ there is a unique vector field
$\xi$, called the {\bf Reeb vector field}, satisfying the two
conditions
$$\xi \hook \eta=1, \qquad \xi \hook d\eta =0.$$
\end{lemma}

\begin{definition}\label{almostcontactdef}
An {\bf almost contact structure} on a differentiable manifolds
$M$ is a triple $(\xi,\eta,\Phi),$ where $\Phi$ is a tensor field
of type $(1,1)$ ({\it i.e.}, an endomorphism of $TM$), $\xi$ is a
vector field, and $\eta$ is a $1$-form which satisfy
$$\eta(\xi)=1 ~~~\hbox{and}~~~ \Phi\circ \Phi= -\BOne + \xi\otimes\eta,$$
where $\BOne$ is the identity endomorphism on $TM.$ A smooth
manifold with such a structure is called an {\bf almost contact
manifold}.
\end{definition}

\begin{remark}\label{contactalmostconremark}
 The reader will notice from Definitions \ref{contactdef} and \ref{almostcontactdef} that an almost contact structure actually has more structure than a contact structure! This is in stark contrast to the usual relationship between a structure and its `almost structure'; however, we feel that the terminology is too well ensconced in the literature to be changed at this late stage. 
\end{remark}

Let $(M,\eta)$ be a contact manifold with a contact $1$-form
$\eta$ and consider $\cald = \ker~\eta\subset TM.$  The subbundle
$\cald$ is maximally non-integrable and it is called the {\it
contact distribution}. The pair $(\cald,\gro)$, where $\gro$ is
the restriction of $d\eta$ to $\cald$ gives $\cald$ the structure
of a symplectic vector bundle. We denote by $\calj(\cald)$ the
space of all almost complex structures $J$ on $\cald$ that are
compatible with $\gro,$ that is the subspace of smooth sections
$J$ of the endomorphism bundle ${\rm End}(\cald)$ that satisfy
\begin{equation}
J^2= -\BOne, \ \ d\eta(JX,JY)=d\eta(X,Y),\ \ d\eta(JX,X)>0
\end{equation}
for any smooth sections $X,Y$ of $\cald.$ Notice that each $J\in
\calj(\cald)$ defines a Riemannian metric $g_\cald$ on $\cald$ by
setting
\begin{equation}\label{tranmetric}
g_\cald(X,Y) =d\eta(JX,Y).
\end{equation}
One easily checks that $g_\cald$ satisfies the compatibility
condition $g_\cald(JX,JY)=g_\cald(X,Y).$ Furthermore, the map
$J\mapsto g_\cald$ is one-to-one, and the space $\calj(\cald)$ is
contractible. A choice of $J$ gives $M$ an almost CR structure.
Moreover, by extending $J$ to all of $TM$ one obtains an almost
contact structure. There are some choices of conventions to make
here. We define the section $\Phi$ of ${\rm End}(TM)$ by $\Phi =J$
on $\cald$ and $\Phi\xi=0$, where $\xi$ is the Reeb vector field
associated to $\eta.$ We can also extend the transverse metric
$g_\cald$ to a metric $g$ on all of $M$ by
\begin{equation}\label{sametric}
g(X,Y)= g_\cald +\eta(X)\eta(Y)= d\eta(\Phi X,Y)+ \eta(X)\eta(Y)
\end{equation}
for all vector fields $X,Y$ on $M.$ One easily sees that $g$
satisfies the compatibility condition $g(\Phi X,\Phi
Y)=g(X,Y)-\eta(X)\eta(Y).$

\begin{definition}\index{contact metric structure}
A contact manifold $M$ with a contact form $\eta$, a vector field
$\xi,$ a section $\Phi$ of ${\rm End}(TM),$ and a Riemannian
metric $g$ which satisfy the conditions
$$\eta(\xi)=1,\qquad \Phi^2=-\BOne +\xi\otimes \eta,$$
$$g(\Phi X,\Phi Y)
=g(X,Y)-\eta(X)\eta(Y)$$ is known as a {\bf metric contact
structure} on $M.$
\end{definition}

\begin{defn-thm}\label{sas.def}
\index{Sasakian manifold}\index{Sasakian structure} A Riemannian
manifold $(M,g)$ is called a {\bf Sasakian manifold} if any one,
hence all, of the following equivalent conditions hold:
\begin{enumerate}
\item[(i)] There exists a Killing vector field $\xi$ of unit
length on $M$ so that the tensor field $\Phi$ of type $(1,1)$,
defined by $\Phi(X) ~=~ -\nabla_X \xi$, satisfies the condition
$$(\nabla_X \Phi)(Y) ~=~ g(X,Y)\xi -g(\xi,Y)X$$
for any pair of vector fields $X$ and $Y$ on $M.$ \item[(ii)]
There exists a Killing vector field $\xi$ of unit length on $M$ so
that the Riemann curvature satisfies the condition
$$R(X,\xi)Y ~=~ g(\xi,Y)X-g(X,Y)\xi,$$
for any pair of vector fields $X$ and $Y$ on $M.$ \item[(iii)] The
metric cone $(\calc(M), \bar{g}) ~=~ (\bbr_+\times M, \
dr^2+r^2g)$ is K\"ahler.
\end{enumerate}
\end{defn-thm}

We refer to the quadruple $\cals=(\xi,\eta,\Phi,g)$ as a {\it
Sasakian structure} on $M$, where $\eta$ is the 1-form dual vector
field $\xi.$ It is easy to see that $\eta$ is a contact form whose
Reeb vector field is $\xi$. In particular
$\cals=(\xi,\eta,\Phi,g)$ is a special type of {\it metric contact
structure}.

\index{characteristic foliation} The vector field $\xi$ is nowhere
vanishing, so there is a 1-dimensional foliation $\calf_\xi$
associated with every Sasakian structure, called the {\it
characteristic foliation}. We will denote the space of leaves of
this foliation by $\calz$. Each leaf of  $\calf_\xi$ has a
holonomy group associated to it.  The dimension of the closure of
the leaves is called the {\it rank} of $\cals$.  We shall be
interested in the case $\hbox{rk}(\cals)=1$. We have

\begin{definition}\label{quasi-regular}\index{Sasakian structure!regular}
The characteristic foliation $\calf_\xi$ is said to be {\bf
quasi-regular}\index{Sasakian structure!quasi-regular} if there is
a positive integer $k$ such that each point has a foliated
coordinate chart $(U,x)$ such that each leaf of $\calf_\xi$ passes
through $U$ at most $k$ times. Otherwise $\calf_\xi$ is called
{\bf irregular}\index{Sasakian structure!irregular}. If $k=1$ then
the foliation is called {\bf regular}\index{Sasakian
structure!regular}, and we use the terminology {\bf non-regular}
to mean quasi-regular, but not regular.
\end{definition}

Let $(M,\cals)$ be a Sasakian manifold, and consider the contact
subbundle $\cald=\ker~\eta.$ There is an orthogonal splitting of
the tangent bundle as
\begin{equation}
TM=\cald \oplus L_\xi,
\end{equation}
where $L_\xi$ is the trivial line bundle generated by the Reeb
vector field $\xi.$ The contact subbundle $\cald$ is just the
normal bundle to the characteristic foliation $\calf_\xi$
generated by $\xi.$ It is naturally endowed with both a complex
structure $J=\Phi|_\cald$ and a symplectic structure $d\eta.$
Hence, $(\cald,J,d\eta)$ gives $M$ a {\it transverse K\"ahler}
structure with K\"ahler form $d\eta$ and metric $g_\cald$ defined
as in (\ref{tranmetric}) which is related to the Sasakian metric
$g$ by $g=g_\cald \oplus \eta\otimes \eta$ as in (\ref{sametric}).
We have  \cite{BG00a} the following fundamental structure theorem:

\begin{theorem}
Let $(M,\xi,\eta,\Phi,g)$ be a compact quasi-regular Sasakian
manifold of dimension $2n+1$, and let $\calz$ denote the space of
leaves of the characteristic foliation. Then the leaf space
$\calz$ is a Hodge orbifold with K\"ahler metric $h$ and K\"ahler
form $\gro$ which defines an integral class $[\gro]$ in
$H^2_{orb}(\calz,\bbz)$ so that $\pi:(M,g) \ra{1.3} (\calz,h)$ is
an orbifold Riemannian submersion. The fibers of $\pi$ are totally
geodesic submanifolds of $M$ diffeomorphic to $S^1.$
\end{theorem}

\noindent and its converse:

\begin{theorem} Let $(\calz,h)$ be a Hodge orbifold.
Let $\pi:M\ra{1.3} \calz$ be the $S^1$ V-bundle whose first Chern
class is $[\gro],$ and let $\eta$ be a connection $1$-form in $M$
whose curvature is $2\pi^*\gro,$ then $M$ with the metric
$\pi^*h+\eta\otimes\eta$ is a Sasakian orbifold.  Furthermore, if
all the local uniformizing groups inject into the group of the
bundle $S^1,$ the total space $M$ is a smooth Sasakian manifold.
\end{theorem}

Irregular structures can be understood by the following result of
Rukimbira \cite{Ruk95a}:

\begin{theorem}\label{Rukapproxthm}
Let $(\xi,\eta,\Phi,g)$ be a compact irregular Sasakian structure
on a manifold $M.$ Then the group $\gA\gu\gt(\xi,\eta,\Phi,g)$ of
Sasakian automorphisms contains a torus $T^k$ of dimension $k\geq
2.$ Furthermore, there exists a sequence
$(\xi_i,\eta_i,\Phi_i,g_i)$ of quasi-regular Sasakian structures
that converge to $(\xi,\eta,\Phi,g)$ in the $C^\infty$
compact-open topology.
\end{theorem}

The orbifold cohomology groups $H_{orb}^p(\calz,\bbz)$ were
defined by Haefliger \cite{Hae84}. In analogy with the smooth case
a {\it Hodge orbifold} is then defined to be a compact K\"ahler
orbifold whose K\"ahler class lies in $H_{orb}^2(\calz,\bbz)$.
Alternatively, we can develop the concept of basic cohomology
which works equally well in the irregular case, but only has
coefficients in $\bbr.$ It is nevertheless quite useful in trying
to extend the notion of $\calz$ being Fano to both the
quasi-regular and the irregular situation. This can be done in
several ways. Here we will use the notion of basic Chern classes.
Recall \cite{Ton} that a smooth p-form $\gra$ on $M$ is called
{\it basic} if
\begin{equation}
\xi\hook \gra=0, \qquad \call_\xi\gra=0\,,
\end{equation}
and we let $\grL^p_B$ denote the sheaf of germs of basic $p$-forms
on $M,$ and by $\grO_B^p$ the set of global sections of $\grL^p_B$
on $M.$ The sheaf $\grL^p_B$ is a module over the ring,
$\grL^0_B,$ of germs of smooth basic functions on $M.$ We let
$C^\infty_B(M)=\grO^0_B$ denote global sections of $\grL^0_B,$
i.e. the ring of smooth basic functions on $M.$  Since exterior
differentiation preserves basic forms we get a de Rham complex
\begin{equation}
\cdots\ra{2.5}\grO_B^p\fract{d}{\ra{2.5}}\grO_B^{p+1}\ra{2.5}\cdots
\end{equation}
whose cohomology $H^*_B(\calf_\xi)$ is called the {\it basic
cohomology} of $(M,\calf_\xi).$ The basic cohomology ring
$H^*_B(\calf_\xi)$ is an invariant of the foliation $\calf_\xi$
and hence, of the Sasakian structure on $M.$ It is related to the
ordinary de Rham cohomology $H^*(M,\bbr)$ by the long exact
sequence \cite{Ton}
\begin{equation}\label{exact}
\cdots\ra{2.5}H_B^p(\calf_\xi)\ra{2.5}H^p(M,\bbr)\fract{j_p}{\ra{2.5}}
H_B^{p-1}(\calf_\xi) \fract{\grd}{\ra{2.5}}
H^{p+1}_B(\calf_\xi)\ra{2.5}\cdots
\end{equation}
where $\grd$ is the connecting homomorphism given by
$\grd[\gra]=[d\eta\wedge \gra]=[d\eta]\cup[\gra],$ and $j_p$ is
the composition of the map induced by $\xi\hook$ with the
well-known isomorphism $H^r(M,\bbr)\approx H^r(M,\bbr)^{S^1}$
where $H^r(M,\bbr)^{S^1}$ is the $S^1$-invariant cohomology
defined from the  $S^1$-invariant $r$-forms $\grO^r(M)^{S^1}.$ We
also note that $d\eta$ is basic even though $\eta$ is not. Next we
exploit the fact that the transverse geometry is K\"ahler. Let
$\cald_\bbc$ denote the complexification of $\cald,$ and decompose
it into its eigenspaces with respect to $J,$ that is, $\cald_\bbc=
\cald^{1,0}\oplus \cald^{0,1}.$ Similarly, we get a splitting of
the complexification of the sheaf $\grL^1_B$ of basic one forms on
$M,$ namely
$$\grL^1_B\otimes \bbc = \grL^{1,0}_B\oplus \grL^{0,1}_B\,.$$
We let $\cale^{p,q}_B$ denote the sheaf of germs of basic forms of
type $(p,q),$ and we obtain a splitting
\begin{equation}
\grL^r_B\otimes \bbc = \bigoplus_{p+q=r}\cale^{p,q}_B\,.
\end{equation}

The basic cohomology groups $H^{p,q}_B(\calf_\xi)$ are fundamental
invariants of a Sasakian structure which enjoy many of the same
properties as the ordinary Dolbeault cohomology of a K\"ahler
structure.

Consider the complex vector bundle $\cald$ on a Sasakian manifold
$(M,\xi,\eta,\Phi,g).$ As such $\cald$ has Chern classes
$c_1(\cald),\ldots,c_n(\cald)$ which can be computed by choosing a
connection $\nabla^\cald$ in $\cald$ \cite{Kobbook}. Let us choose
a local  foliate unitary transverse frame $(X_1,\ldots,X_n),$ and
denote by $\grO^T$ the transverse curvature $2$-form with respect
to this frame. A simple calculation shows that $\grO^T$ is a basic
$(1,1)$-form. Since the curvature $2$-form $\grO^T$ has type
$(1,1)$ it follows as in ordinary Chern-Weil theory that

\begin{theorem}\label{transChernWeil}
The $k^{\rm th}$ Chern class $c_k(\cald)$ of the complex vector
bundle $\cald$ is represented by the basic $(k,k)$-form $\grg_k$
determined by the formula
$$\det\Bigl(\BOne_n-\frac{1}{2\pi i}\grO^T\Bigr)=1+\grg_1+\cdots +\grg_k.$$
Since $\grg_k$ is a closed basic $(k,k)$-form it represents an
element in $H_B^{k,k}(\calf_\xi)\subset H_B^{2k}(\calf_\xi)$ that
is called the {\it basic $k^{\rm th}$ Chern class} and denoted by
$c_k(\calf_\xi).$
\end{theorem}

We now concentrate on the first Chern classes $c_1(\cald)$ and
$c_1(\calf_\xi).$ We have

\begin{definition}\label{c1def}\index{Sasakian structure!null}
\index{Sasakian structure!positive} \index{Sasakian
structure!negative}A Sasakian structure $\cals=(\xi,\eta,\Phi,g)$
is said to be {\bf positive (negative)} if $c_1(\calf_\xi)$ is
represented by a positive (negative) definite $(1,1)$-form. If
either of these two conditions is satisfied $\cals$ is said to be
{\bf definite}, and otherwise $\cals$ is called {\bf indefinite}.
$\cals$ is said to be {\bf null} if $c_1(\calf_\xi)=0.$
\end{definition}

Notice that irregular structures cannot occur for negative or null
Sasakian structures, since the dimension of
$\gA\gu\gt(\xi,\eta,\Phi,g)$ is greater than one. In analogy with
common terminology of smooth algebraic varieties we see that a
positive Sasakian structure is a {\it transverse Fano
structure}\footnote{For a more algebro-geometric approach to
positivity and fundamentals on log Fano orbifolds  see
\cite{BG05}.}, while a null Sasakian structure is a {\it
transverse Calabi-Yau structure}. The negative Sasakian case
corresponds to the canonical bundle being ample.

\medskip
\section{Sasaki-Einstein and 3-Sasakian Geometry}\label{se.3se}
\medskip\index{Sasaki-Einstein manifold}

\begin{definition}
A Sasakian manifold $(M,\cals)$ is {\bf Sasaki-Einstein} if the
metric $g$ is also Einstein. 
\end{definition}

For any 2n+1-dimensional Sasakian
manifold ${\rm Ric}(X,\xi)=2n\eta(X)$ implying that any
Sasaki-Einstein metric must have positive scalar curvature. Thus
any complete Sasaki-Einstein manifold must have a finite
fundamental group. Furthermore the metric cone
$(\calc(M),\bar{g})= (\bbr_+\times M, \ dr^2+r^2g)$ on $M$ is
K\"ahler Ricci-flat (Calabi-Yau).

The following theorem \cite{BG00a} is an orbifold version of the
famous Kobayashi bundle construction of Einstein metrics on
bundles over positive K\"ahler-Einstein manifolds \cite{Bes,
Kob56}.

\begin{theorem}\label{kob}
 Let $(\calz, h)$ be a compact Fano orbifold
with $\pi_1^{orb}(\calz)=0$ and K\"ahler-Einstein metric $h$. Let
$\pi:M\ra{1.3} \calz$ be the $S^1$ V-bundle whose first Chern
class is $\frac{c_1(\calz)}{\hbox{Ind}(\calz)}.$ Suppose further
that the local uniformizing groups of $\calz$ inject into $S^1.$
Then with the metric $g=\pi^*h+\eta\otimes\eta$, $M$ is a compact
simply connected Sasaki-Einstein manifold.
\end{theorem}

Here ${\hbox{Ind}(\calz)}$ is the {\it orbifold Fano index}
\cite{BG00a} defined to be the largest positive integer such that
$\frac{c_1(\calz)}{\hbox{Ind}(\calz)}$ defines a class in the
orbifold cohomology group $H^2_{orb}(\calz,\bbz).$  A very special
class of Sasaki-Einstein spaces is naturally related to several
quaternionic geometries.

\begin{definition}\label{3s.def}\index{$3$-Sasakian manifold}
Let $(M,g)$ be a Riemannian manifold of dimension $m$. We say that
$(M,g)$ is {\bf $3$-Sasakian} if the metric cone $(\calc(M),\bar{g})=
(\bbr_+\times M, \ dr^2+r^2g)$ on $M$  is hyperk\"ahler.
\end{definition}

We emphasize the important observation of Kashiwada \cite{Kas71} that a 3-Sasakian manifold is automatically Einstein.
We denote a Sasakian manifold with a $3$-Sasakian structure
by $(M,\boldsymbol{\cals}),$ where
$\boldsymbol{\cals}=(\cals_1,\cals_2,\cals_3)$ is a triple or a
$2$-sphere of Sasakian structures
$\cals_i=(\eta_i,\xi_i,\Phi_i,g).$

\begin{remark} In the $3$-Sasakian case there is an extra structure, {\it i.e.},
the transverse geometry $\calo$ of the $3$-dimensional foliation
which is quaternionic-K\"ahler. In this case, the transverse space
$\calz$ is the twistor space of $\calo$ and the natural map
$\calz\ \ra{1.3}\ \calo$ is the orbifold twistor fibration
\cite{Sal82}. We get the following diagram which we denote by
$\boldsymbol{\diamondsuit}(M,\boldsymbol{\cals})$ \cite{BGM93a,
BGM94a}: \vskip 20pt
\begin{center}
\begin{tabular}{|c|}\hline {\bf Hyperk\"ahler}\\{\bf Geometry}\\ \hline\end{tabular}
\end{center}
\begin{equation}
\begin{tabular}{|c|}\hline {\bf Twistor}\\{\bf Geometry}\\ \hline
\end{tabular}
\begin{array}{ccccc}
&&\calc(M)&&\\
&\swarrow&&\searrow&\\
\calz &&\hskip -15pt\la{4}\hskip -30pt\decdnar{}&& M\\
&\searrow& &\swarrow&\\
&&\calo&&
\end{array}
\begin{tabular}{|c|}\hline {\bf 3-Sasakian}\\{\bf Geometry}\\ \hline
\end{tabular}
\end{equation}
\begin{center}
\begin{tabular}{|c|}\hline {\bf Quaternion K\"ahler}\\{\bf Geometry}\\ \hline\end{tabular}
\end{center}
\end{remark}

\begin{remark} The table below summarizes properties of the cone and transverse
geometries associated to various metric contact structures. \vskip
12pt
\begin{center}
  \begin{tabular}{|c|c|c|}
  \hline
    Cone Geometry of $\calc(M)$ & $M$ & Transverse Geometry of $\calf_\xi$ \\ \hline\hline
    Symplectic & Contact & Symplectic \\ \hline
    K\"ahler &Sasakian & K\"ahler \\ \hline
    K\"ahler &positive Sasakian & Fano, $c_1(\calz)>0$ \\ \hline
    K\"ahler &null Sasakian & Calabi-Yau, $c_1(\calz)=0$ \\ \hline
    K\"ahler &negative Sasakian & ample canonical bundle, $c_1(\calz)<0$ \\ \hline
    Calabi-Yau &Sasaki-Einstein & Fano, K\"ahler-Einstein \\ \hline
    Hyperk\"ahler &3-Sasakian& $\bbc$-contact, Fano, K\"ahler-Einstein \\ \hline
  \end{tabular}
\end{center}
\end{remark}
\bigskip
For numerous examples and constructions of \Se and $3$-Sasakian
manifolds see \cite{BG05}. We finish this section with a remark
that both the $3$-Sasakian metric on $M$ and the twistor space
metric on $\calz$ admit `squashings' which are again Einstein.
More generally, let $\pi:M\longrightarrow B$ be an orbifold
Riemannian submersion with $g$ the Riemannian metric on $M.$ Let
$\calv$ and $\calh$ denote the vertical and horizontal subbundles
of the tangent bundle $TM.$ For each real number $t>0$ we
construct a one parameter family $g_t$ of Riemannian metrics on
$M$ by defining
\begin{equation} g_t|_{\calv} =tg|_{\calv}\,, \qquad g_t|_\calh =g|_\calh\,, \qquad
g_t(\calv,\calh)=0\,.
\end{equation}
So for each $t>0$ we have an orbifold Riemannian submersion with
the same base space. Furthermore, if the fibers of $g$ are totally
geodesic, so are the fibers of $g_t\,.$ We apply the canonical
variation to the orbifold Riemannian submersion
$\pi:M\longrightarrow \calo$ and  $\pi:\calz\longrightarrow \calo$

\begin{theorem} \label{secondE.thm}
Every $3$-Sasakian manifold $M$ admits a second Einstein metric of
positive scalar curvature. Furthermore, the twistor space $\calz$
also admits a second orbifold Einstein metric which is
Hermitian-Einstein, but not K\"ahler-Einstein.
\end{theorem}

\medskip
\section{Toric Sasaki-Einstein $5$-Manifolds}\label{toric.SE.mfds}
\medskip\index{Sasaki-Einstein manifold!toric}
Examples of \Se manifolds are plentiful and we refer the
interested reader to our monograph for a detailed exposition
\cite{BG05}. Here we would like to consider the toric \Se
structures in dimension $5$ again referring to \cite{BG05} for all
necessary details. Toric \Se $5$-manifolds recently emerged from
physics in the context of {\it supersymmetry}  and the
so-called AdS/CFT duality conjecture which we will discuss in the
last section. It is known that, in dimension $5$, toric \Se
structures can only occur on the $k$-fold connected sums
$k(S^2\times S^3)$ \cite{BG05}. The first inhomogeneous toric \Se
structures on $S^2\times S^3$ were constructed by Gauntlett,
Martelli, Sparks, and Waldram. It follows that $S^2\times S^3$
admits infinitely many distinct quasi-regular and irregular toric
\Se structures \cite{GMSW04a}. Toric geometry of these examples
was further explored in \cite{MaSp05b, MaSpYau05, MaSpYau06}. We
will now describe a slightly different approach to a more general
problem.

Consider the symplectic reduction of $\bbc^n$ (or equivalently the
Sasakian reduction of $S^{2n-1}$) by a $k$-dimensional torus
$T^k.$ Every complex representation of a $T^k$ on $\bbc^{n}$ can
be described by an exact sequence
$$0\ra{1.3}T^k\fract{f_\Omega}{\ra{1.3}} T^{n}\ra{1.3} T^{n-k}\ra{1.3}0\, .$$
The monomorphism $f_\grO$ can be represented by the diagonal
matrix
$$f_\Omega(\tau_1,\ldots,\tau_k)={\rm diag}\Biggl(
\displaystyle{\prod_{i=1}^k\tau_i^{a^i_1}},\ldots,
\displaystyle{\prod_{i=1}^k\tau_i^{a^i_{n}}}\Biggr)\, ,$$ where
$(\tau_1,..,\tau_k)\in S^1\times\cdots\times S^1=T^k$ are the
complex coordinates on $T^k,$ and $a^i_\gra\in \bbz$  are the
coefficients of a $k\times n$ integral {\it weight} matrix
$\Omega\in\calm_{k,n}(\bbz).$ We have \cite{BG05}

\begin{proposition}\label{toric.Sasakian.quotients}
Let $X(\Omega)=(\bbc^n\setminus{\bf0})/\!\!/T^k(\Omega)$ denote
the K\"ahler quotient of the standard flat K\"ahler structure on
$(\bbc^n\setminus{\bf0})$  by the weighted Hamiltonian
$T^k$-action with an integer weight matrix $\Omega.$ Consider the
K\"ahler moment map
\begin{equation}
\mu^i_\Omega(\bfz)=
\sum_{\alpha=1}^{n}a^i_{\alpha}|z_\alpha|^2,\qquad i=1,\ldots,k\,.
\end{equation}
If all minor $k\times k$ determinants of \ $\Omega$ are non-zero
then $X(\Omega)=C(Y(\Omega))$ is a cone on a compact Sasakian
orbifold $Y(\Omega)$ of dimension $2(n-k)-1$ which is the Sasakian
reduction of the standard Sasakian structure on $S^{2n-1}.$ In
addition, the projectivization of $X(\Omega)$ defined by
$\calz(\Omega)=X(\Omega)/\bbc^*$ is a K\"ahler reduction of the
complex projective  space $\bbc\bbp^{n-1}$ by a Hamiltonian
$T^k$-action defined by $\Omega$ and it is the transverse space of
the Sasakian structure on $Y(\Omega)$ induced by the quotient. If
\begin{equation}\label{CY.condition}
\sum_\alpha a^i_\alpha=0,\qquad \forall~ i=1,\ldots,k
\end{equation}
then $c_1(X(\Omega))=c_1(\cald)=0.$ In particular, the orbibundle
$Y(\Omega)\ra{1.2}\calz(\Omega)$ is anticanonical. Moreover, the
cone $C(Y(\Omega))$, its Sasakian base $Y(\Omega)$, and the
transverse space $\calz(\Omega)$ are all toric orbifolds.
\end{proposition}

\begin{remark}\label{reductproprem}
The conditions on the matrix $\grO$ that assure that $Y(\grO)$ is
a smooth manifold are straightforward to work out. They involve
gcd conditions on certain minor determinants of $\grO.$
\end{remark}

This proposition is nicely summarized by the `reduction' diagram
\begin{equation}\label{commquot4}
\begin{matrix} \bbc\bbp^{n-1} &\longleftarrow  &S^{2n-1}  &\longleftarrow & \bbc^n \setminus (\bf0) \\
                      \Downarrow  &   &\Downarrow &         &\Downarrow \\
                  \calz(\grO) &\longleftarrow   &Y(\grO) &\longleftarrow
                  &C(Y(\grO)).
\end{matrix}
\end{equation}

Both the toric geometry and the topology of $Y(\Omega)$ depend on
$\Omega$. Furthermore, $Y(\Omega)$ comes equipped with a family of
Sasakian structures. When $n-k=3,$ assuming that $Y(\Omega)$ is
simply connected (which is an additional condition on $\Omega$),
we must have $m(S^2\times S^3)$ for some $m\leq k.$ We will be
mostly interested in the case when $m=k.$

Gauntlett, Martelli, Sparks, and Waldram \cite{GMSW04a} gave an
explicit construction of a \Se metric for $\grO=(p,p,-p+q,-p-q),$
where $p$ and $q$ are relatively prime nonnegative integers with
$p>q.$ (The general case for $k=1$ was treated later in
\cite{CLPP05, MaSp05b}, see Remark \ref{CLPP.metrics} below). To
connect with the original notation we write $Y(\Omega)=Y_{p,q}.$
Then we get:

One can check that  $Y_{1,0}$ is just the homogeneous metric on
$S^2\times S^3$ which is both toric and regular. The next simplest
example is $Y_{2,1}$ which, as a toric contact (Sasakian)
manifold, is a circle bundle over the blow up of $\bbc\bbp^2$ at
one point $F_1=\bbc\bbp^2\#\overline{\bbc\bbp^2}$ \cite{MaSp06}.
As $F_1$ cannot admit any K\"ahler-Einstein metric, Kobayashi's
bundle construction cannot give a compatible Sasaki-Einstein
structure. But there is a choice of a Reeb vector field in the
torus which makes it possible to give $Y_{2,1}$ a \Se metric. The
\Se structure on $Y_{2,1}$ is not quasi-regular and this was the
first such example in the literature. Hence, $S^2\times S^3$
admits infinitely many  toric quasi-regular \Se structures and
infinitely many toric irregular \Se structures of rank $2.$ We
have the following generalization of the $Y_{p,q}$ metrics due to
\cite{FOW06, CFO07}:

\begin{theorem}\label{toric.SE.exist.unique} Let $Y(\Omega)$ be as in Proposition
\ref{toric.Sasakian.quotients}. Then $Y(\Omega)$ admits a toric
\Se structure which is unique up to a transverse biholomorphism.
\end{theorem}

This existence of a \Se metric is proved in \cite{FOW06} although
the authors do not draw all the conclusions regarding possible
toric \Se manifolds that can be obtained. They give one
interesting example of an irregular \Se structure which
generalizes the $Y_{2,1}$ example of \cite{MaSp05b} in the
following sense: One considers a regular positive Sasakian
structure on the anticanonical circle bundle over the del Pezzo
surface $\bbc\bbp^2\#2\overline{\bbc\bbp^2}$ which gives a toric
Sasakian structure on $2(S^2\times S^3).$ The regular Sasakian
structure on $2(S^2\times S^3)$ cannot have any \Se metric.
However, as it is with $Y_{2,1}$ Futaki, Ono and Wang \cite{FOW06}
show that one can deform the regular structure to a unique
irregular \Se structure. A slightly different version of the
Theorem \ref{toric.SE.exist.unique} is proved in \cite{CFO07}
where uniqueness is also established. Cho, Futaki and Ono work
with toric diagrams rather than with K\"ahler (Sasakian) quotients
which amounts to the same thing by Delzant's construction. We
should add that the results of \cite{CFO07} apply to the toric \Se
manifolds in general dimension and not just in dimension $5$.

\begin{corollary}\label{toric.SE.5man} The manifolds
$k(S^2\times S^3)$ admit infinite families of toric  \Se
structures for each $k\geq1.$
\end{corollary}

As in the $k=1$ case one would expect infinitely many
quasi-regular and infinitely many irregular such \Se structures
for each $\Omega$ satisfying all the condition.

\begin{remark}\label{CLPP.metrics}
The general anticanonical circle reduction was considered
independently in two recent papers, \cite{CLPP05,MaSp05b}. There
it was shown that for $\Omega=\bfp=(p_1,p_2,-q_1,-q_2)$, with
$p_i,q_i\in\bbz^+$, $p_1+p_2=q_1+q_2$, and ${\rm gcd}(p_i,q_j)=1$
for all $i,j=1,2,$ the $5$-manifold $Y(\Omega)\approx S^2\times
S^3$ admits a \Se structure which coincides with that on $Y_{p,q}$
when $p_1=p_2=p$ and $q_1=p-q,q_2=p+q.$  In \cite{CLPP05} this
family is denoted by $L^5(a,b,c),$ where $\bfp=(a,b,-c,-a-b+c)$
and they write the metric explicitly. However, in this case it
appears to be harder (though, in principle, possible) to write
down the condition under which the \Se Reeb vector field
$\xi=\xi(a,b,c)$ is quasi-regular. A priori, it is not even clear
whether the quasi-regularity condition has any additional
solutions beyond those obtained for the subfamily $Y_{p,q}.$
Moreover, it follows from \cite{CFO07} that the metrics of
\cite{CLPP05,MaSp05b} describe all possible toric \Se structures
on $S^2\times S^3.$

There have been similar constructions of a two-parameter family
$X_{p,q}$ of toric \Se metrics on $2(S^2\times S^3)$ \cite{HKW05},
and another two-parameter family, called $Z_{p,q},$ on
$3(S^2\times S^3)$ \cite{OoYa06}. All these examples, and many
more, can be obtained as special cases of Theorem
\ref{toric.SE.exist.unique} as they are all $Y(\Omega)$ for some
choice of $\Omega$. The $Y_{p,q},$ $L^5(a,b,c),$ $X_{p,q}$ and
$Z_{p,q}$ metrics have received a lot of attention because of the
role such \Se manifolds play in the AdS/CFT Duality Conjecture.
They created an avalanche of papers studying the properties of
these metrics from the physics perspective
\cite{ABCC06,OoYa06,OoYa06c,OoYa06b,
KSY05,HEK05,BuZa05,BBC05,BFZ05,SfZo05,HKW05,BLMP05,
BHK05,BFHMS05,Pal05, HSY04}. The AdS/CFT duality will be discussed
in the last section.
\end{remark}

\medskip
\section{The Dirac Operator and Killing Spinors}\label{Diracsect}
\medskip

We begin with a definition of spinor bundles and the bundle of
Clifford algebras of a vector bundle \cite{LaMi89,Fri00}. Recall
that the {\it Clifford algebra} $Cl(\bbr^n)$ over $\bbr^n$ can be
defined as the quotient algebra of the tensor algebra
$\calt(\bbr^n)$ by the two-sided ideal $\cali$ generated by
elements of the form $v\otimes v +q(v)$ where $q$ is a quadratic
form on $\bbr^n.$

\begin{definition}\label{Clifford}
Let  $E$ be a vector bundle with inner product $\langle\cdot,\cdot
\rangle$ on a smooth manifold $M$, and let $\calt(E)$ denote the
tensor bundle over $E.$  The {\bf Clifford bundle} of $E$ is the
quotient bundle $Cl(E)=\calt(E)/\cali(E)$ where $\cali$ is the
bundle of ideals (two-sided) generated pointwise by elements of
the form $v\otimes v+\langle v,v\rangle$ with $v\in E_x\,.$ A {\bf
real spinor bundle} $S(E)$ of $E$ is a bundle of modules over the
Clifford bundle $Cl(E).$ Similarly, a {\bf complex spinor bundle}
is a bundle of complex modules over the complexification
$Cl(E)\otimes \bbc.$
\end{definition}

As vector bundles $Cl(E)$ is isomorphic to the exterior bundle
$\grL(E),$ but their algebraic structures are different. The
importance of $Cl(E)$ is that it contains the spin group
$Spin(n)$, the universal (double) covering group of the orthogonal
group $SO(n),$ so one obtains all the representations of $Spin(n)$
by studying representations of  $Cl(E).$ We assume that the vector
bundle $E$ admits a spin structure, so $w_2(E)=0.$ We are
interested mainly in the case when $(M,g)$ is a Riemannian spin
manifold and $E=TM$ in which case we write $S(M)$ instead of
$S(TM).$ The Levi-Civita connection $\nabla$ on $TM$ induces a
connection, also denoted $\nabla,$ on any of the spinor bundles
$S(M),$ or more appropriately on the sections $\grG(S(M)).$

\begin{definition}\label{Dirac}
Let $(M^n,g)$ be a Riemannian spin manifold and let $S(M)$ be any
spinor bundle. The {\bf Dirac operator} is the first order
differential operator $D:\grG(S(M))\ra{1.8} \grG(S(M))$ defined by
$$D\psi=\sum_{j=1}^n E_j\cdot \nabla_{E_j}\psi\, ,$$
where $\{E_j\}$ is a local orthonormal frame and $\cdot$ denotes
Clifford multiplication.
\end{definition}

The Dirac operator, of course originating with the famous Dirac
equation describing fermions  in theoretical physics, was brought
into mathematics by Atiyah and Singer in \cite{AtSi63}. Then
Lichnerowicz \cite{Lic63b} proved his famous result that a
Riemannian spin manifold with positive scalar curvature must have
vanishing $\hat{A}$-genus. An interesting question on any spin
manifold is: what are the eigenvectors of the Dirac operator. In
this regard the main objects of interest consists of special
sections of certain spinor bundles called {\it Killing spinor
fields} or just {\it Killing spinors} for short. Specifically,
(cf. \cite{BFGK91,Fri00})

\begin{definition}\label{real.Killing.spinors}\index{Killing
spinor!real} \index{Killing spinor!parallel} \index{Killing
spinor!imaginary}Let $(M,g)$ be a complete $n$-dimensional
Riemannian spin manifold, and let $S(M)$ be a spin bundle (real or
complex) on $M$ and $\psi$ a smooth section of $S(M).$ We say that
$\psi$ is a {\bf Killing spinor} if for every vector field $X$
there is $\gra\in \bbc,$ called {\bf Killing number}, such that
$$\nabla_X\psi=\alpha X\!\cdot\!\psi\, .$$
Here $X\!\cdot\!\psi$ denotes the Clifford product of $X$ and
$\psi.$ We say that $\psi$ is {\bf imaginary} when $\alpha\in {\rm
Im}(\bbc^*),$ $\psi$ is {\bf parallel} if $\alpha=0$ and $\psi$ is
{\bf real}\footnote{Here the standard terminology real and
imaginary Killing spinors can be somewhat misleading. The Killing
spinor $\psi$ is usually a section of a complex spinor bundle. So
a real Killing spinor just means that $\gra$ is real.} if
$\alpha\in {\rm Re}(\bbc^*).$\index{Killing spinor}
\end{definition}

We shall see shortly that the three possibilities for the Killing
number $\gra:$ real, imaginary, or $0,$ are the only
possibilities. The name Killing spinor derives from the fact that
if $\psi$ is a non-trivial Killing spinor and $\gra$ is real, the
vector field
\begin{equation}\label{KspinKvector}
X_\psi =\sum_{j=1}^n g(\psi,E_j\cdot \psi)E_j
\end{equation}
is a Killing vector field for the metric $g$ (which, of course,
can be zero). If $\psi$ is a Killing spinor on an $n$-dimensional
spin manifold, then
\begin{equation}\label{KspineigenD}
D\psi= \sum_{j=1}^n E_j\cdot \nabla_{E_j}\psi =\sum_{j=1}^n\alpha
E_j\cdot E_j\cdot\!\psi= -n\gra \psi\, .
\end{equation}
So Killing spinors are eigenvectors of the Dirac operator with
eigenvalue $-n\gra.$ In 1980 Friedrich \cite{Fr80} proved the
following remarkable theorem:

\begin{theorem}\label{thm.fried}
Let $(M^n,g)$ be a Riemannian spin manifold which admits a
non-trivial Killing spinor $\psi$ with Killing number $\gra.$ Then $(M^n,g)$ is Einstein with
scalar curvature $s=4n(n-1)\gra^2.$
\end{theorem}

A proof of this is a straightforward curvature computation which
can be found in either of the books \cite{BFGK91,Fri00}. It also
uses the fact that a non-trivial Killing spinor vanishes nowhere.
It follows immediately from Theorem \ref{thm.fried} that $\gra$
must be one of the three types mentioned in Definition
\ref{real.Killing.spinors}. So if the Killing number is real then
$(M,g)$ must be a positive Einstein manifold. In particular, if
$M$ is complete, then it is compact. On the other hand if the
Killing number is pure imaginary, Friedrich shows that $M$ must be
non-compact.

The existence of Killing spinors not only puts restrictions on the
Ricci curvature, but also on both the Riemannian and the Weyl
curvature operators \cite{BFGK91}.

\begin{proposition}\label{K.S.curvature} Let $(M^n,g)$ be a Riemannian
spin manifold. Let $\psi$ be a Killing spinor on $M$ with Killing
number $\alpha$ and let $\calr,
\calw:\Lambda^2M\ra{1.2}\Lambda^2M$ be the Riemann and Weyl
curvature operators, respectively. Then for any vector field $X$
and any $2$-form $\beta$ we have
\begin{align} &\calw(\beta)\cdot\psi=0\, ; \\
&(\nabla_X\calw)(\beta)\cdot\psi=-2\alpha\bigl(X\hook\calw(\beta)\bigr)\cdot\psi\, ; \\
&(\calr(\beta)+4\alpha^2\beta)\cdot\psi=0\, ;\\
&(\nabla_X\calr)(\beta)\cdot\psi=-2\alpha\bigl(X\hook\calr(\beta)+4\alpha^2\beta(X)\bigr)
\cdot\psi\, .
\end{align}
\end{proposition}

These curvature equations can be used to prove (see \cite{BFGK91}
or \cite{Fri00})

\begin{theorem}\label{K.S.irreducibility}
Let $(M^n,g)$ be a connected Riemannian spin manifold admitting a
non-trivial Killing spinor with $\alpha\not=0.$ Then $(M,g)$ is
locally irreducible. Furthermore, if $M$ is locally symmetric, or
$n\leq4,$ then $M$ is a space of constant sectional curvature
equal to $4\alpha^2.$
\end{theorem}

Friedrich's main objective in \cite{Fr80} was an improvement of
Lichnerowicz's estimate in \cite{Lic63b} for the eigenvalues of
the Dirac operator. Indeed, Friedrich proves that the eigenvalues
$\grl$ of the Dirac operator on any compact manifold satisfy the
estimate
\begin{equation}\label{fried.bound}
\lambda^2\geq \frac{1}{4}\frac{ns_0}{n-1}\, ,
\end{equation}
where $s_0$ is the minimum of the scalar curvature on $M.$ Thus,
Killing spinors $\psi$ are eigenvectors that realize equality in
equation \eqref{fried.bound}.  Friedrich also proves the converse
that any eigenvector of $D$ realizing the equality must be a
Killing spinor with
\begin{equation}\label{mineigen}
\gra =\pm \frac{1}{2}\sqrt{\frac{s_0}{n(n-1)}}\, .
\end{equation}

\begin{example}\label{Kspinsphere}{\bf [Spheres]}
In the case of the round sphere $(S^n,g_0)$ equality in equation
\eqref{fried.bound} is always attained. So normalizing such that
$s_0=n(n-1),$ and using B\"ar's Correspondence Theorem
\ref{Barcorresthm} below the number of corresponding real Killing
spinors equals the number of constant spinors on $\bbr^{n+1}$ with
the flat metric. The latter is well known (see the appendix of
\cite{PeRi88}) to be $2^{\lfloor n/2\rfloor}$ for each of the
values $\gra=\pm\frac{1}{2},$ where $\lfloor n/2\rfloor$ is the
largest integer less than or equal to $n/2.$
\end{example}

\begin{remark} Actually (without making the connection to
Sasakian geometry) already in \cite{Fr80} Friedrich gives a
non-spherical example of a compact $5$-manifold with a real
Killing spinor: $M=SO(4)/SO(2)$ with its homogeneous
Kobayashi-Tanno \Se structure.
\end{remark}

We now wish to relate Killing spinors to the main theme of this
article, Sasakian geometry. First notice that if a Sasakian
manifold $M^{2n+1}$ admits a Killing spinor, Theorem
\ref{thm.fried} says it must be Sasaki-Einstein, so the scalar
curvature $s_0=2n(2n+1)$, and equation \eqref{mineigen} implies
that $\gra=\pm \frac{1}{2}.$ We have the following result of
Friedrich and Kath \cite{FrKat2}

\begin{theorem}\label{FrKat.thm}
Every simply connected \Se manifold admits non-trivial real
Killing spinors. Furthermore,
\begin{enumerate}
\item[(i)] if $M$ has dimension $4m+1$ then $(M,g)$ admits exactly
one Killing spinor for each of the values $\gra=\pm \frac{1}{2},$
\item[(ii)] if $M$ has dimension $4m+3$ then $(M,g)$ admits at
least two Killing spinors for one of the values $\gra=\pm
\frac{1}{2}.$
\end{enumerate}
\end{theorem}

\begin{proof}[Outline of Proof]
(Details can be found in \cite{FrKat2} or the book \cite{BFGK91}.)
Every simply connected \Se manifold is known to be spin, so $M$
has a spin bundle $S(M).$ Given a fixed Sasakian structure
$\cals=(\xi,\eta,\Phi,g)$ we consider two subbundles
$\cale_\pm(\cals)$ of $S(M)$ defined by
\begin{equation}\label{Rebb.Killing}
\cale_\pm(\cals)=\{\psi\in S(M)\ \ |\ \ (\pm2\Phi X+\pounds_\xi
X)\cdot\psi=0,\ \ \ \forall X\in\Gamma(TM)\}\, .
\end{equation}
Set $\nabla^\pm_X= \nabla_X \pm \frac{1}{2}X\cdot.$ A
straightforward computation shows that $\nabla^\pm$ preserves the
subbundles $\cale_\pm$ and defines a connection there. Moreover,
by standard curvature computations it can be shown that the
connection $\nabla^\pm$ is flat in $\cale_\pm(\cals).$ So it has
covariantly constant sections which are precisely the Killing
spinors. One then uses some representation theory of ${\rm
Spin}(2n+1)$ to compute the dimensions of $\cale_+(\cals)$ and
$\cale_-(\cals)$ proving the result.
\end{proof}

We have the following:

\begin{corollary}\label{SE.locally.irr}
Let $(M,g)$ be a \Se manifold of dimension $2m+1.$ Then
$(M,g)$ is locally symmetric if and only if $(M,g)$ is of constant
curvature. Moreover, ${\rm Hol}(g)=SO(2m+1)$ and $(M,g)$ is
locally irreducible as a Riemannian manifold.
\end{corollary}

\begin{proof} If necessary, go to the universal cover
$\tilde M$. This is a compact simply connected \Se manifold;
hence, it admits a non-trivial Killing spinor by Theorem
\ref{FrKat.thm}. The first statement then follows from the Theorem
\ref{K.S.irreducibility}. The second statement follows from the
Berger Theorem \ref{Berger.main.thm}. Since $M$ has dimension
$2m+1$ the only possibilities for ${\rm Hol}(g)$ are $SO(2m+1)$
and $G_2$. But the latter is Ricci flat, so it cannot be
Sasaki-Einstein.
\end{proof}

Friedrich and Kath began their investigation in dimension $5$
\cite{FrKa89} where they showed that a simply-connected compact
$5$-manifold which admits a Killing spinor must be
Sasaki-Einstein. In dimension $7$ they showed that there are
exactly three possibilities: weak $G_2$-manifolds, \Se manifolds
which are not $3$-Sasakian, and $3$-Sasakian manifolds
\cite{FrKat2}.  Later Grunewald gave a description of
$6$-manifolds admitting Killing spinors \cite{Gru90}. We should
add an earlier result of Hijazi who showed that the only
$8$-dimensional manifold with Killing spinors must be the round
sphere \cite{Hij86}. By 1990 a decade of research by many people
slowly identified all the ingredients of a classification of such
manifolds in terms of their underlying geometric structures. The
pieces of the puzzle consisting of round spheres in any dimension,
\Se manifolds in odd dimensions, nearly K\"ahler manifolds in
dimension $6$, and weak $G_2$-holonomy manifolds in dimension $7$
were all in place with plenty of interesting examples to go around
\cite{BFGK91}. What remained at that stage was to show that in
even dimensions greater than $8$ there is nothing else but the
round spheres, while in odd dimensions greater than $7$ the only
such examples must be Sasaki-Einstein. The missing piece of the
puzzle was finally uncovered by B\"ar:  real Killing spinors on
$M$ correspond to parallel spinors on the cone $C(M)$
\cite{Bar93}. A bit earlier Wang \cite{Wan89} had shown that on a
simply connected complete Riemannian spin manifold the existence
of parallel spinors corresponds to reduced holonomy.  This led
B\"ar to an elegant description of the geometry of manifolds
admitting real Killing spinors (in any dimension) in terms of
special holonomies of the associated cones. We refer to the
correspondence between real Killing spinors on $M$ and parallel
spinors on the cone $C(M)$ (equivalently reduced holonomy) as {\it
B\"ar's correspondence}. In particular, this correspondence not
only answered the last remaining open questions, but also allowed
for simple unified proofs of most of the theorems obtained
earlier.

\section{Real Killing Spinors, Holonomy and B\"ar's
Correspondence}\index{B\"ar's correspondence} As mentioned the
B\"ar correspondence relates real Killing spinors on a compact
Riemannian spin manifold $(M,g)$ to parallel spinors on the
Riemannian cone $(C(M),\bar{g}).$ We now make this statement
precise.
\begin{theorem}\label{Barcorresthm}
Let $(M^n,g)$ be a complete Riemannian spin manifold and
$(C(M^n),\bar{g})$ be its Riemannian cone. Then there is a one to
one correspondence between real Killing spinors on $(M^n,g)$ with
$\gra=\pm\frac{1}{2}$ and parallel spinors on $(C(M^n),\bar{g}).$
\end{theorem}

\begin{proof}
The existence of a parallel spinor on $(C(M^n),\bar{g})$ implies
that $\bar{g}$ is Ricci flat by Theorem \ref{thm.fried}. Then by
Lemma \ref{einlemma} $(M^n,g)$ is Einstein with scalar curvature
$s=n(n-1).$ So any Killing spinors must have $\gra=\pm\frac{1}{2}$
by equation \eqref{mineigen}. As in the proof of Theorem
\ref{FrKat.thm}, $\nabla^\pm_X= \nabla_X \pm \frac{1}{2}X\cdot$
defines a connection in the spin bundle $S(M).$ The connection
$1$-forms $\gro^\pm$ of $\nabla^\pm$ are related to the connection
$1$-form  $\gro$ of the Levi-Civita connection by $\gro^\pm =\gro
\pm \frac{1}{2}\grb,$ where $\grb$ is a $1$-form called the {\it
soldering form}. This can be interpreted as a connection with
values in the Lie algebra $\gs\gp\gi\gn(n+1)=
\gs\gp\gi\gn(n)\oplus \bbr^n,$ and pulls back to the Levi-Civita
connection in the spin bundles on the cone $(C(M^n),\bar{g})$. So
parallel spinors on the cone correspond to parallel spinors on
$(M,g)$ with respect to the connection $\nabla^\pm$ which
correspond precisely to real Killing spinors with respect to the
Levi-Civita connection.
\end{proof}

Now we have the following definition:

\begin{definition}\label{p.q.type}\index{Killing spinor!type}
We say that a Riemannian spin manifold $(M,g)$ is of {\bf type
$(p,q)$} if it carries exactly $p$ linearly independent real
Killing spinors with $\alpha>0$ and exactly $q$ linearly
independent real Killing spinors with $\alpha<0$.
\end{definition}

The following theorem has an interesting history. As mentioned
above it was B\"ar \cite{Bar93} who recognized the correspondence
between real Killing spinors on $(M,g)$ and parallel spinors on
the Riemannian cone $(C(M),\bar{g})$. The relation between
parallel spinors and reduced holonomy was anticipated in the work
of Hitchin \cite{Hit74a} and Bonan \cite{Bon66}, but was
formalized in the 1989 paper of Wang \cite{Wan89}. It has also
been generalized to the non-simply connected case in
\cite{Wan95,MoSe00}.

\begin{theorem}\label{Bar.main.thm1}
Let $(M^n,g)$ be a complete simply connected Riemannian spin
manifold, and let ${\rm Hol}(\bar{g})$ be the holonomy group of
the Riemannian cone $(C(M),\bar{g})$. Then $(M^n,g)$ admits a
non-trivial real Killing spinor with $(M^n,g)$ of type $(p,q)$ if
and only if $\bigl(\dim M, {\rm Hol}(\bar{g}), (p,q)\bigr)$ is one
of the $6$ possible triples listed in the table below:

\begin{center}\vbox{\[\begin{array}{|l|l|l|} \hline \dim (M)&{\rm Hol}(\bar{g})
&{\rm type}~(p,q) \hl\hline n&{\rm id}&(2^{\lfloor
n/2\rfloor},2^{\lfloor n/2\rfloor})\hl 4m+1&SU(2m+1)&(1,1)\hl
4m+3&SU(2m+2)&(2,0)\hl 4m+3&Sp(m+1)&(m+2,0)\hl 7& Spin(7)&(1,0)\hl
6&G_2&(1,1)\hl
\end{array}\]}
\end{center}
\noindent Here $m\geq 1,$ and $n>1.$
\end{theorem}

\begin{proof}[Outline of Proof]
Since $(M,g)$ is complete and has a non-trivial real Killing
spinor, it is compact by Theorem \ref{thm.fried}. It then follows
from a theorem of Gallot \cite{Gal79} that if the Riemannian cone
$(C(M),\bar{g})$ has reducible holonomy it must be flat. So we can
apply Berger's Theorem \ref{Berger.main.thm}. Now Wang
\cite{Wan89} used the spinor representations of the possible
irreducible holonomy groups on Berger's list to give the
correspondence between these holonomy groups and the existence of
parallel spinors. First he showed that the groups listed in Table
\ref{bergertable} that are not on the above table do not admit
parallel spinors. Then upon decomposing the spin representation of
the group in question into irreducible pieces, the number of
parallel spinors corresponds to the multiplicity of the trivial
representation. Wang computes this in all but the first line of
the table when $(C(M),\bar{g})$ is flat. In this case $(M,g)$ is a
round sphere as discussed in Example \ref{Kspinsphere}, so the
number of linearly independent constant spinors is $(2^{\lfloor
n/2\rfloor},2^{\lfloor n/2\rfloor}).$ By B\"ar's Correspondence
Theorem \ref{Barcorresthm} real Killing spinors on $(M,g)$
correspond precisely to parallel spinors on $(C(M),\bar{g}).$ Note
that the hypothesis of completeness in Wang's theorem \cite{Wan89}
is not necessary, so that the correspondence between the holonomy
groups and parallel spinors holds equally well on Riemannian
cones. However, the completeness assumption on $(M,g)$ guarantees
the irreducibility of the cone $(C(M),\bar{g})$ as mentioned
above.
\end{proof}

Let us briefly discuss the types of geometry involved in each case
of this theorem. As mentioned in the above proof the first line of
the table corresponds to the round spheres. The next three lines
correspond to \Se geometry, so Theorem \ref{Bar.main.thm1}
generalizes the Friedrich-Kath Theorem \ref{FrKat.thm} in this
case. The last of these three lines corresponds precisely to
3-Sasakian geometry by Definition \ref{3s.def}. Finally the two
cases whose cones have exceptional holonomy will be discussed in
more detail in Section \ref{weakG2sect} below. Suffice it here to
mention that it was observed by Bryant and Salamon \cite{BrSa89}
that a cone on a nearly parallel $G_2$ manifold has its own holonomy
in $Spin(7).$ It is interesting to note that Theorem
\ref{Bar.main.thm1} generalizes the result of Hijazi in dimension
eight mentioned earlier as well as part of the last statement in
Theorem \ref{K.S.irreducibility}, namely

\begin{corollary}\label{Bar.main.thm0}
Let $(M^{2n},g)$ be a complete simply connected Riemannian spin
manifold of dimension $2n$ with $n\neq 3$ admitting a non-trivial
real Killing spinor. Then $M$ is isometric to the round sphere.
\end{corollary}

We end this section with a brief discussion of the non-simply
connected case. Here we consider two additional cases for ${\rm
Hol}(\bar{g})$, namely $SU(2m+2)\rtimes\bbz_2$ and $Sp(2)\times
\bbz_d.$ See \cite{Wan95,MoSe00} for the list of possibilities.

\begin{example}\label{Zmod2Stiefel}
${\rm Hol}(\bar{g})=SU(2m)\rtimes \bbz_2$. Consider the
$(4m-1)$-dimensional Stiefel manifold $V_2(\bbr^{2m+1})$ with its
homogeneous \Se metric. The quotient manifold $M^{4m-1}_\grs$ of
$V_2(\bbr^{2m+1})$ by the free involution $\grs$ induced from
complex conjugation has an Einstein metric which is ``locally
Sasakian''. The cone $C(M^{4m-1}_\grs)$ is not K\"ahler and its
holonomy is ${\rm Hol}(\bar{g})=SU(2m+2)\rtimes\bbz_2$. According
to Wang \cite{Wan95} $C(M^{4m-1}_\grs)$ admits a spin structure
with precisely one parallel spinor if and only if $m$ is even, and
according to Moroianu and Semmelmann \cite{MoSe00}
$C(M^{4m-1}_\grs)$ admits exactly two spin structures each with
precisely one parallel spinor if $m$ is even. Thus, by Theorem
\ref{Barcorresthm} $M^{4m-1}_\grs$ admits exactly two spin
structures each with exactly one Killing spinor if and only if $m$
is even.
\end{example}

\begin{example}\label{3Sas->SE}
Consider a $3$-Sasakian manifold $(M^{4n-1},\boldsymbol{\cals})$
and choose a Reeb vector field $\xi(\boldsymbol{\grt}).$ Let $C_m$
be the cyclic subgroup of order $m>2$ of the circle group
generated by $\xi(\boldsymbol{\grt}).$ Assume that $m$ is
relatively prime to the order $\upsilon(\boldsymbol{\cals})$ of
$\boldsymbol{\cals}$ and that the generic fibre of the fundamental
$3$-dimensional foliation $\calf_Q$ is $SO(3),$ so that $C_m$ acts
freely on $M^{4n-1}.$ This last condition on the generic fibre is
easy to satisfy; for example, it holds for any of the $3$-Sasakian
homogeneous spaces other than the standard round sphere, as well
as the bi-quotients described in \cite{BGM94a}. (To handle the
case when the generic fibre is $Sp(1)$ we simply need to divide
$m$ by two when it is even). Since $C_m$ is not in the center of
$SO(3),$ the quotient $M^{4n-1}/C_m$ is not $3$-Sasakian. However,
$C_m$ does preserve the Sasakian structure determined by
$\xi(\boldsymbol{\grt}),$ so $M^{4n-1}/C_m$ is Sasaki-Einstein.
The cone $C(M^{4n-1}/C_m)$ has holonomy $Sp(n)\times \bbz_m$, and
admits precisely $\frac{n+1}{m}$ parallel spinors if and only if
$m$ divides $n+1$ \cite{Wan95,MoSe00}. Thus, by Theorem
\ref{Barcorresthm} $M^{4n-1}/C_m$ admits precisely $\frac{n+1}{m}$
Killing spinors when  $m$ divides $n+1.$
\end{example}

\medskip
\section{Geometries Associated with 3-Sasakian 7-manifolds}\label{lbspin3}
\medskip

It is most remarkable that to each $4n$-dimensional positive QK
metric $(\calo,g_\calo)$ (even just locally) one can associate
{\it nine} other Einstein metrics in dimensions $4n+k$,
$k=1,2,3,4$. Alternatively, one could say that each $3$-Sasakian
metric $(M,g)$ canonically defines an additional nine Einstein
metrics in various dimensions. We have already encountered all of
these metrics. First there are the four geometries of the diamond
diagram $\boldsymbol{\diamondsuit}(M,\boldsymbol{\cals}).$ Then
$M$ and $\calz$ admit additional ``squashed" Einstein metrics
discussed in Theorem \ref{secondE.thm}. Thus we get five Einstein
metrics with positive Einstein constants: $(\calo,g_\calo), (M,g),
(M',g'), (\calz, h), (\calz', h')$. Of course $M\simeq M'$ and
$\calz\simeq \calz'$ as smooth manifolds (orbifolds) but they are
different as Riemannian manifolds (orbifolds), hence, the
notation. Let us scale all these metrics so that the Einstein
constant equals the dimension of the total space minus 1. Note
that any $3$-Sasakian metric already has this property. In the
other four cases this is a choice of scale which is quite natural
due to Lemma \ref{einlemma}. However, note that this is not the
scale one gets for $(\calz,h),$ and $(\calo,g_\calo)$ via the
Riemannian submersion from $(M,g).$ Now, in each case one can
consider its Riemannian cone which will be Ricci-flat by Lemma
\ref{einlemma}.
%$(C(\calo),
%dt^2+t^2g_\calo), (C(M), dt^2+t^2g), (C(M), dt^2+t^2g'),
%(C(\calz), dt^2+t^2h), (C(\calz'), dt^2+t^2h')$
We thus obtain five Ricci-flat metrics on the corresponding
Riemannian cones. In addition, one can also take (iterated)
sine-cone metrics defined in \eqref{conemetric} on the same five
bases. These metrics are all Einstein of positive scalar curvature
(cf. Lemma \ref{sine-coneEinstein}). Let us summarize all this
with the following extension of
$\boldsymbol{\diamondsuit}(M,\boldsymbol{\cals})$:

\begin{equation}\label{10.Einstein.metrics}
\xymatrix{&C(\calz')\ &\ar@{_{(}->}@<.2ex>[l]\ \calz'\ar[rd]&&M'\ar[ld]\ \ar@{^{(}->}@<-.2ex>[r]&\
C(M')\\
&&&\calo\ \ar@{^{(}->}@<-.2ex>[r]&\ \ C(\calo)\\
&C(\calz)\ &\ar@{_{(}->}@<.2ex>[l]\ \calz\ar[ru]&&\ar[ll]M\ar[lu]\
\ar@{^{(}->}@<-.2ex>[r]&\ C(M)}\end{equation}

There would perhaps be nothing special about all these $10$ (and
many more by iterating sine-cone construction) geometries beyond
what has already been discussed in the previous sections. This is
indeed true when ${\rm dim}(M)>7.$ However, when ${\rm dim}(M)=7$,
or, alternatively, when $\calo$ is a positive self-dual Einstein
orbifold metric (more generally, just a local metric of this type)
some of the metrics occurring in diagram
\eqref{10.Einstein.metrics} have additional properties. We shall
list all of them first. For the moment, let us assume that $(M,g)$
is a compact $3$-Sasakian $7$-manifold, then the following hold:
\begin{enumerate}
\item $(\calo,g_\calo)$ is a positive self-dual Einstein manifold
(orbifold). We will think of it as the {\it source} of all the
other geometries. \item $(C(\calo),dt^2+t^2g_\calo)$ is a
$5$-dimensional Ricci-flat cone with  base $\calo.$ \item $(\calz,
h)$ is the orbifold twistor space of $\calo.$ \item $(\calz', h')$
is a nearly-K\"ahler manifold (orbifold). \item $(M,g)$ is the
$3$-Sasakian manifold. \item $(M',g')$ is a $7$-manifold with weak
$G_2$ structure.\item $(C(\calz'),dt^2+t^2h')$ is a $7$-manifold
with holonomy inside $G_2.$ \item $(C_s(\calz'),dt^2+(\sin
^2t)h')$ is a $7$-manifold with weak $G_2$ structure. \item
$(C(\calz),dt^2+t^2h)$ is a $7$-dimensional Ricci-flat cone with
base $\calz.$  \item $(C(M),dt^2+t^2g)$ is hyperk\"ahler with
holonomy contained in $Sp(2).$ \item $(C(M'),dt^2+t^2g')$ has
holonomy contained in $Spin(7).$
\end{enumerate}
The cases (2) and (8) do not appear to have any special properties
other than Ricci-flatness. The cases (1), (3), (5), and (10) are
the four geometries of
$\boldsymbol{\diamondsuit}(M,\boldsymbol{\cals}).$ The five
remaining cases are all very interesting from the point of view of
the classification of Theorem \ref{Bar.main.thm1}. Indeed $\calz'$
and $C(\calz')$ are examples of the structures listed in the last
row of the table while $C_2(\calz')$, $M'$ and $C(M')$ give
examples of the structures listed in the fifth row. In particular,
our diagram \eqref{10.Einstein.metrics} provides for a cornucopia
of the orbifold examples in the first case and smooth manifolds in
the latter.

\subsection{Nearly Parallel $G_2$-Structures and $Spin(7)$ Holonomy
Cones}\label{weakG2sect}\index{nearly parallel $G_2$-structures} Recall, that
geometrically $G_2$ is defined to be the Lie group acting on the
imaginary octonions $\bbr^7$ and preserving the $3$-form
\begin{equation}
\begin{split}
\varphi&=\alpha_1\wedge\alpha_2\wedge\alpha_3+
\alpha_1\wedge(\alpha_4\wedge\alpha_5- \alpha_6\wedge\alpha_7) \\
&+ \alpha_2\wedge(\alpha_4\wedge\alpha_6- \alpha_7\wedge\alpha_5)+
\alpha_3\wedge(\alpha_4\wedge\alpha_7-
\alpha_5\wedge\alpha_6),\label{bg2.2.1}
\end{split}
\end{equation}
where $\{\alpha_i\}_{i=1}^7$ is a fixed orthonormal basis of the
dual of $\bbr^7$. A $G_2$ structure on a $7$-manifold $M$ is, by
definition, a reduction of the structure group of the tangent
bundle to $G_2$. This is equivalent to the existence of a global
3-form $\varphi\in\Omega^3(M)$ which may be written locally as
\ref{bg2.2.1}.  Such a 3-form defines an associated Riemannian
metric, an orientation class, and a spinor field of constant
length. 

\begin{definition}\label{weak.G2}
Let $(M,g)$ be a complete $7$-dimensional Riemannian manifold. We
say that $(M,g)$ is a {\bf nearly parallel}\footnote{It had become customary to refer to this notion as `weak holonomy $G_2$', a terminology introduced by Gray \cite{Gra71}. However, it was pointed out to us by the anonomous referee that this terminology is misleading due to the fact that Gray's paper contains errors rendering the concept of weak holonomy useless as discovered by Alexandrov \cite{Ale05}. Hence, the term `nearly parallel' used in \cite{FKMS97} is preferred.} $G_2$ structure if there exist a
global $3$-form $\varphi\in\Omega^3(M)$ which locally can be
written in terms of a local orthonormal basis as in \ref{bg2.2.1},
and $d\varphi=c\star \varphi$, where $\star$ is the Hodge star
operator associated to $g$ and $c\neq 0$ is a constant whose sign is
fixed by an orientation convention.  \end{definition}

The case $c=0$ in Definition \ref{weak.G2} is somewhat
special. In particular, it is known \cite{Sal89} that the
condition $d\varphi =0=d\star\varphi$ is equivalent to the
condition that $\varphi$ be parallel, {\it i.e.},
$\nabla\varphi=0$ which is equivalent to the condition that the
metric $g$ has holonomy group contained in $G_2.$  The following
theorem provides the connection with the previous discussion on
Killing spinors \cite{Bar93}

\begin{theorem}\label{weak.G2.Spin7}
Let $(M,g)$ be a complete $7$-dimensional Riemannian manifold with
a nearly parallel $G_2$ structure. Then the holonomy ${\rm Hol}(\bar{g})$ of the
metric cone $(C(M),\bar{g})$ is contained in $Spin(7).$ In
particular, $C(M)$ is Ricci-flat and $M$ is Einstein with positive
Einstein constant $\lambda=6.$
\end{theorem}

\begin{remark}
The sphere $S^7$ with its constant curvature metric is isometric
to the isotropy irreducible space $Spin(7)/G_2.$ The fact that
$G_2$ leaves invariant (up to constants) a unique $3$-form and a
unique $4$-form on $\bbr^7$ implies immediately that this space
has a nearly parallel $G_2$ structure.
\end{remark}

\begin{definition}\label{proper.weak.G2}
 Let $(M,g)$ be a complete $7$-dimensional Riemannian manifold.
We say that $g$ is a {\bf proper $G_2$-metric} if ${\rm
Hol}(\bar{g})={\rm Spin}(7).$
\end{definition}

We emphasize here that $G_2$ is the structure group of $M,$ not
the Riemannian holonomy group. Specializing Theorem
\ref{Bar.main.thm1} to dimension $7$ gives the following theorem
due to Friedrich and Kath \cite{FrKat2}.

\begin{theorem}\label{killing.spinor.7D}
Let $(M^7,g)$ be a complete simply-connected Riemannian spin
manifold of dimension $7$ admitting a non-trivial real Killing
spinor with $\alpha>0$ or $\alpha<0$. Then there are four
possibilities:
\begin{enumerate}
\item[(i)] $(M^7,g)$ is of type $(1,0)$ and it is a proper
$G_2$-manifold, \item[(ii)] $(M^7,g)$ is of type $(2,0)$ and it is
a Sasaki-Einstein manifold, but $(M^7,g)$ is not $3$-Sasakian,
\item[(iii)] $(M^7,g)$ is of type $(3,0)$ and it is $3$-Sasakian,
\item[(iv)] $(M^7,g)=(S^7,g_{can})$ and is of type $(8,8).$
\end{enumerate}
Conversely, if $(M^7,g)$ is a compact simply-connected proper
$G_2$-manifold then it carries precisely one Killing spinor with
$\alpha>0.$ If $(M^7,g)$ is a compact simply-connected \Se
$7$-manifold which is not $3$-Sasakian then $M$ carries precisely
$2$ linearly independent Killing spinors with $\alpha>0$. Finally,
if $(M^7,g)$ is a $3$-Sasakian $7$-manifold, which is not of
constant curvature, then $M$ carries precisely $3$ linearly
independent Killing spinors with $\alpha>0.$
\end{theorem}

\begin{remark}\label{complete.G2}
The four possibilities of the Theorem \ref{killing.spinor.7D}
correspond to the sequence of inclusions
$${\rm Spin}(7) \supset SU(4)\supset Sp(2)\supset \BOne\, .$$ 
All of the
corresponding cases are examples of nearly parallel $G_2$ metrics.
If we exclude the trivial case when the associated cone is flat,
we have three types of nearly parallel $G_2$ geometries. Following
\cite{FKMS97} we use the number of linearly independent Killing
spinors to classify these geometries, and 
call them type I, II, and III corresponding to cases (i),
(ii), and (iii) of Theorem \ref{killing.spinor.7D}, respectively.
\end{remark}

We are now ready to describe the $G_2$ geometry of the
$M'\hookrightarrow C(M')$ part of Diagram
\ref{10.Einstein.metrics} \cite{GaSal96, FKMS97}:

\begin{theorem}\label{Strick.G2.from.3Sas}
Let $(M,\boldsymbol{\cals})$ be a $7$-dimensional $3$-Sasakian
manifold. Then the $3$-Sasakian metric $g$ is a nearly parallel
$G_2$ metric. Moreover, the second Einstein metric $g'$ given by Theorem
\ref{secondE.thm} and scaled so that the Einstein constant
$\lambda=6$ is a nearly parallel $G_2$ metric; in fact, it is a proper
$G_2$ metric.
\end{theorem}

\begin{proof} For the second Einstein metric $g'$ we have three
mutually orthonormal $1$-forms
$\alpha^1=\sqrt{t}\eta^1,\quad\alpha^2=\sqrt{t}\eta^2,\quad
\alpha^3=\sqrt{t}\eta^3,$ where $t$ is the parameter of the
canonical variation. Let $\{\alpha^4,\alpha^5,\alpha^6,\alpha^7\}$
be local $1$-forms spanning the annihilator of the vertical
subbundle $\calv_3$ in $T^*\cals$ such that
$$\hi^1= 2(\alpha^4\wedge\alpha^5-\alpha^6\wedge\alpha^7)\, ,$$
$$\hi^2= 2(\alpha^4\wedge\alpha^6-\alpha^7\wedge\alpha^5)\, ,$$
$$\hi^3= 2(\alpha^4\wedge\alpha^7-\alpha^5\wedge\alpha^6)\, .$$
Then the set $\{\gra^1,\ldots, \gra^7\}$ forms a local orthonormal
coframe for the metric $g'.$ Let
\begin{equation}
\Eta = \eta_1\wedge\eta_2\wedge \eta_3\,,\qquad
\Theta=\sum_a\eta_a\wedge \hi_a=\sum_a\eta_a\wedge d\eta_a+
6\Eta\label{btop.2.2}
\end{equation}
In terms of the $3$-forms $\Eta$ and $\Theta$ we have
$\varphi=\frac{1}{2}\sqrt{t}\Theta+\sqrt{t}^3\Eta.$ One easily
sees that this is of the type of equation \eqref{bg2.2.1} and,
therefore, defines a compatible $G_2$-structure.  Moreover, a
straightforward computation gives
$$d\varphi =\frac{1}{2}\sqrt{t}\Omega + \sqrt{t}(t+1)d\Eta,\qquad \star\varphi =
-\frac{1}{2}td\Eta - \frac{1}{24}\Omega\, .$$ Thus,
$d\varphi=c\star\varphi$ is solved with $\sqrt{t}=1/\sqrt5$, and
$c=-12/\sqrt5$. So $g'$ is nearly parallel. That $g'$ is a
proper $G_2$ metric is due to \cite{FKMS97}. The idea is to use
Theorem \ref{killing.spinor.7D}. Looking at the four possibilities
given in that theorem, we see that it suffices to show that $g'$
is not Sasaki-Einstein. The details are in \cite{FKMS97}.
\end{proof}

\begin{example} $3$-Sasakian $7$-manifolds
are plentiful \cite{BG05}. All of them give, by Theorem
\ref{Strick.G2.from.3Sas}, examples of type I and type III
geometries. Examples of simply connected type I geometries that do
not arise via Theorem \ref{Strick.G2.from.3Sas} are the
homogeneous Aloff-Wallach spaces $M^7_{m,n},$ $(m,n)\not=(1,1)$ which,
as special cases of Eschenburg bi-quotients \cite{CMS96, BFGK91}, are
together with an isotropy irreducible homogeneous space defined as
follows: Consider the space $\calh_2$ of homogeneous polynomials
of degree $2$ in three real variables $(x_1,x_2,x_3).$ As ${\rm
dim}(\calh_2)=5$ it gives rise to the embedding $SO(3)\subset
SO(5).$ We take $M=SO(5)/SO(3).$ This example was used by Bryant
to get the first $8$-dimensional metric with holonomy $Spin(7)$
\cite{Bry87}. Examples of type II geometries (Sasaki-Einstein) are
equally rich \cite{BG05}. In particular, there are hundreds of
examples of type II nearly parallel $G_2$ metrics on each of the 28
homotopy spheres in dimension $7.$
\end{example}
\begin{remark}
According to \cite{CMS96} the Aloff-Wallach manifold $M^7_{1,1}$
has three Einstein metrics. One is the homogeneous $3$-Sasakian
metric. The second is the proper $G_2$ metric of Theorem
\ref{Strick.G2.from.3Sas}. The third Einstein metric is also
nearly parallel most likely being of type I, but we could not
positively exclude type II as a possibility.
\end{remark}

\begin{oproblem} Classify all compact $7$-manifolds with nearly parallel
$G_2$ structures of type I, II, or III, respectively.
\end{oproblem}
The classification of type III consists of the classification of all compact
$3$-Sasakian $7$-manifolds. This is probably very hard. The case
of $3$-Sasakian $7$-manifolds with vanishing
$\ga\gu\gt(M,\boldsymbol{\cals})$ appears quite difficult. The type II
classification ($7$-dimensional \Se manifolds which are not
$3$-Sasakian) is clearly completely out of reach at the moment. A
classification of proper nearly parallel $G_2$ structures on a compact
manifold that do not arise via Theorem \ref{Strick.G2.from.3Sas}
would be very interesting and it is not clear how hard this
problem really is.

\begin{remark}\label{complete.hol.Spin7} The holonomy
$Spin(7)$ cone metrics are plentiful but never complete. However,
some of these metrics can be deformed to complete holonomy
$Spin(7)$ ones on non compact manifolds. The first example was
obtained by Bryant and Salamon who observed that the spin bundle
over $S^4$ with its canonical metric carries a complete metric
with holonomy $Spin(7)$ \cite{BrSa89}. Locally the metric was
later considered also in \cite{GPP90}. More generally, spin
orbibundles over positive QK orbifolds also carry such complete
orbifold metrics as observed by Bryant and Salamon in
\cite{BrSa89}. Other complete examples were constructed later by
physicists \cite{CGLP02, CGLP04, KaYa02a, KaYa02b}. Finally, the
first compact examples were obtained in 1996 by Joyce
\cite{Joy96b, Joy99}. See Joyce's book \cite{Joy00} for an
excellent detailed exposition of the methods and the discussion of
examples.
\end{remark}

\begin{oproblem}\label{3S.complete.Spin7} {\bf [Complete metrics on cones]}
Let $(M^7,\boldsymbol{\cals})$ be any $3$-Sasakian $7$-manifold
and let $(M^7,g')$ be the associated proper nearly parallel $G_2$
squashed metric. Consider the two Riemannian cones for these
metrics.
\begin{enumerate} \item[(i)] When does the metric cone $(C(M),dt^2+t^2g')$ admit  complete
holonomy $Spin(7)$ deformations? \item[(ii)] When does the metric
cone $(C(M),dt^2+t^2g)$ admit  complete holonomy $Sp(2)$
(hyperk\"ahler) deformations?
\end{enumerate}
In other dimensions one also could ask the following more general
questions:
\begin{enumerate}
\item[(iii)] Let $(M^{4n+3},\boldsymbol{\cals})$ be a compact
$3$-Sasakian manifold. When does the metric cone
$(C(M),dt^2+t^2g)$ admit  complete hyperk\"ahler (or just
Calabi-Yau) deformations? \item[(iv)] Let $(M^{2n+1},\cals)$ be a
compact Sasaki-Einstein manifold. When does the metric cone
$(C(M),dt^2+t^2g)$ admit complete Calabi-Yau deformations?
\item[(v)] Let $(M^7,g)$ be a compact nearly parallel
$G_2$-manifold. When does the metric cone $(C(M),dt^2+t^2g)$ admit
complete holonomy $Spin(7)$ deformations? 
\item[(vi)] Let
$(M^6,g)$ be a compact strict nearly K\"ahler manifold. When does
the metric cone $(C(M),dt^2+t^2g)$ admit  complete holonomy $G_2$
deformations?
\end{enumerate}
\end{oproblem}

The metric on the spin bundle $S(S^4)$  by Bryant and Salamon is a
deformation of the $Spin(7)$ holonomy metric on the cone over the
squashed metric on $S^7$ \cite{CGLP02, CGLP04}, so there are
examples of such deformations regarding question (i). Regarding
(ii), we recall that every compact $3$-Sasakian $3$-manifold is
isometric to $S^3/\Gamma$ and the metric cone is the flat cone
$\bbc^2/\Gamma$. Hence, one could think of (ii) as a
$7$-dimensional analogue of a similar problem whose complete
solution was given by Kronheimer \cite{Kro89a}. There are
non-trivial examples also in the higher dimensional cases. The
metric cone on the homogeneous $3$-Sasakian manifold
$\cals(1,1,1)$ of \cite{BGM94a} admits complete hyperk\"ahler
deformations, namely the Calabi metric on $T^*\bbc\bbp^2.$ We do
not know of any other examples at the moment. In case (iv) of
the Calabi-Yau cones on \Se manifolds, however, there are many
such examples. Futaki very recently proved that such a complete Calabi-Yau
metric exists for all the regular toric \Se manifolds of Section
\ref{toric.SE.mfds} \cite{Fut07}. In such cases the metric can be
thought of as a complete Ricci-flat K\"ahler metric on the canonical
bundle over a toric Fano manifold. Futaki's result should
generalize to the case of toric log Fano orbifolds.

\subsection{Nearly K\"ahler $6$-Manifolds and $G_2$ Holonomy
Cones.}\label{nearlyKsect}\index{nearly K\"ahler manifolds} In
this section we explain the geometry of the $\calz'\hookrightarrow
C(\calz')$ part of the diagram \ref{10.Einstein.metrics}. Before
we specialize to dimension 6 we begin with a more general
introduction. Nearly K\"ahler manifolds were first studied by
Tachibana in \cite{Tach59} and they appear under the name of
almost Tachibana spaces in Chapter VIII of the book \cite{Yan65}.
They were then rediscovered by Gray \cite{Gra70} and given the
name nearly K\"ahler manifolds which by now is the accepted name.

\begin{definition}\label{nearly.K}
A {\bf nearly K\"ahler manifold} is  an almost Hermitian manifold
$(M,g,J,\gro)$ such that $(\nabla_XJ)X=0$ for all tangent vectors
$X$, where $\nabla$ is the Levi-Civita connection and $J$ is the
almost complex structure. One says that a nearly K\"ahler manifold
is {\bf strict} if it is not K\"ahler.
\end{definition}

This definition is equivalent to the condition
\begin{equation}\label{JKtensor}
(\nabla_XJ)Y+ (\nabla_YJ)X =0
\end{equation}
for all vector fields $X,Y,$ which is to say that $J$ is a {\it
Killing tensor} field.  An alternative characterization of nearly
K\"ahler manifolds is given by
\begin{proposition}\label{NKcovd}
An almost Hermitian manifold $(M,g,J,\gro)$ is nearly K\"ahler if
and only if
$$\nabla \gro= \frac{1}{3}d\gro\, .$$
In particular, a strict nearly K\"ahler structure is never
integrable.
\end{proposition}

Any nearly K\"ahler manifold can be locally decomposed as the
product of a K\"ahler manifold and a strict nearly K\"ahler
manifold. Such a decomposition is global in the simply connected
case \cite{Nagy02a}. Hence, the study of nearly K\"ahler manifolds
reduces to the case of strict ones. In addition every nearly
K\"ahler manifold in dimension $4$ must be K\"ahler so that the
first interesting dimension is six.

The following theorem establishes relationship between the twistor
space $\calz\ra{1.2}\calo$ of a quaternionic K\"ahler manifold
(orbifold) and nearly K\"ahler geometry.

\begin{theorem}\label{NK.structure.twistor} Let $\pi:(\calz,h)\ra{1.2}(\calo,g_\calo)$
be the twistor space of a positive QK manifold with its K\"ahler
structure $(J,h,\omega_h).$ Then $\calz$ admits a strict nearly
K\"ahler structure $(J_1,h_1,\omega_{h_1})$. If
$TM=\calv\oplus\calh$ is the natural splitting induced by $\pi$
then
\begin{equation}
h\!\!\mid_\calv=2h_1\!\!\mid_\calv\, ,\ \ \
h\!\!\mid_\calh=h_1\!\!\mid_\calh=\pi^*(g_\calo)\, ,\ \ \
\end{equation}
\begin{equation}
J\!\!\mid_\calv=-J_1\!\!\mid_\calv\, ,\ \ \
J\!\!\mid_\calh=J_1\!\!\mid_\calh\, .\ \ \
\end{equation}
\end{theorem}

Theorem \ref{NK.structure.twistor} is due to Eells and Salamon \cite{EeSa83} when
$\calo$ is $4$-dimensional. The higher dimensional analogue was
established in  \cite{AGI98} (see also \cite{Nagy02a}).

\begin{remark} Observe that the metric of the nearly K\"ahler
structure of Theorem \ref{NK.structure.twistor}, in general, is
{\it not} Einstein. In particular, $h_1$ is not the squashed
metric $h'$ introduced in the diagram \ref{10.Einstein.metrics},
unless ${\rm dim}(\calz)=6$. In six dimensions, we can scale $h_1$
so that it has scalar curvature $s=30$ and then indeed $h_1=h'$ as
one can easily check.
\end{remark}

\begin{definition}\label{3symmetric.space} Let $M=G/H$ be a homogeneous space. We say that $M$
is {\bf
$3$-symmetric} if $G$ has an automorphism $\sigma$ of order $3$
such that $G_0^\sigma\subset H\subset G^\sigma$, where $G^\sigma$
is the fixed point set of $\sigma$ and $G_0^\sigma$ is the
identity component in $G_0^\sigma$.
\end{definition}

We have the following two theorems concerning nearly K\"ahler
homogeneous Riemannian manifolds. The first is due to Wolf  and
Gray in all dimensions but six \cite{WoGr68a, WoGr68b}. They also
conjectured that the result is true for strict nearly K\"ahler
$6$-manifolds. The Wolf-Gray conjecture was proved quite recently
by Butruille \cite{But05,But06} which is the second theorem below.

\begin{theorem}\label{Gray.Wolf} Every compact homogeneous
strict nearly K\"ahler manifold $M$ of dimension different than
$6$ is $3$-symmetric.
\end{theorem}

\begin{theorem}\label{But.thm}
Let $(M,g)$ be a strict nearly K\"ahler $6$-dimensional Riemannian
homogeneous manifold. Then $M$ is isomorphic as a homogeneous
space to a finite quotient of $G/H,$ where $G$ and $H$ are one of
the following:
\begin{enumerate}
\item $G=SU(2)\times SU(2)$ and $H=\{{\rm id}\};$ \item $G=G_2$
and $H=SU(3),$ where metrically $G/H=S^6$ the round sphere; \item
$G=Sp(2)$ and $H=SU(2)U(1),$ where $G/H=\bbc\bbp^3$ with its
nearly K\"ahler metric determined by Theorem
\ref{NK.structure.twistor}; \item $G=SU(3)$ and $H=T^2$, where
$G/H$ is the flag manifold with its nearly K\"ahler metric
determined by Theorem \ref{NK.structure.twistor}.
\end{enumerate}

\noindent Each of these manifolds carries a unique invariant
nearly K\"ahler structure, up to homothety.
\end{theorem}

In every dimension, the only known compact examples of nearly
K\"ahler manifolds are $3$-symmetric. On the other hand, Theorem
\ref{NK.structure.twistor} can be easily generalized to the case
of orbifolds so that there are plenty examples of compact
inhomogeneous strict nearly K\"ahler orbifolds in every dimension.

\begin{theorem}\label{nagy.main.class} Let $M$ be a compact
simply-connected strict nearly K\"ahler manifold. Then, in all
dimensions, as a Riemannian manifold $M$ decomposes as a product
of
\begin{enumerate}
\item $3$-symmetric spaces, \item twistor spaces of positive QK
manifolds $\calq$ such that $\calq$ is not symmetric, \item
$6$-dimensional strict nearly K\"ahler manifold other than the
ones listed in Theorem \ref{But.thm}.
\end{enumerate}
\end{theorem}
This theorem is due to Nagy \cite{Nagy02}, but our formulation
uses the result of Butruille together with the fact that the
twistor spaces of all symmetric positive QK manifolds are
$3$-symmetric. The LeBrun-Salamon conjecture can now be phrased as
follows
\begin{conjecture}\label{L.S.NK.conj} Any compact simply connected strict irreducible
nearly K\"ahler manifold $(M,g)$ of dimension greater than $6$
must be a $3$-symmetric space.
\end{conjecture}
In particular, the Conjecture \ref{L.S.NK.conj} is automatically
true in dimensions $4n$ because of Nagy's classification theorem
and also true in dimensions $10$ and $14$ because all positive QK
manifolds in dimension $8$ and $12$ are known. The third case
leads to an important
\begin{oproblem}\label{NK.6.mfds} Classify all compact strict nearly
K\"ahler manifolds in dimension $6$.
\end{oproblem}

Dimension six is special not just because of the r\^ole it plays
in Theorem \ref{nagy.main.class}. They have several remarkable
properties which we summarize in the following theorem.

\begin{theorem}\label{NK6.properties} Let $(M,J,g,\omega_g)$ be a
compact strict nearly K\"ahler $6$-manifold. Then
\begin{enumerate}
\item[(i)] The metric $g$ is Einstein of positive scalar
curvature. \item[(ii)] $c_1(M)=0$ and $w_2(M)=0$. \item[(iii)] If
$g$ is scaled so it has Einstein constant $\lambda=5$ then the
metric cone $(C(M), dt^2t+t^2g)$ has holonomy contained in $G_2$.
In particular, $C(M)$ is Ricci-flat.
\end{enumerate}
\end{theorem}
The first property is due to Matsumoto \cite{Mat72t} while the
second is due to Gray \cite{Gra76}.  The last part is due to B\"ar
\cite{Bar93}. In fact, nearly K\"ahler $6$-manifolds is the
geometry of the last row of the table of Theorem
\ref{Bar.main.thm1}. More precisely we have the following theorem
proved by Grunewald \cite{Gru90}:

\begin{theorem}\label{NK.Gru.thm}
Let $(M^6,g)$ be a complete simply connected Riemannian spin
manifold of dimension $6$ admitting a non-trivial Killing spinor
with $\alpha>0$ or $\alpha<0.$ Then there are two possibilities:
\begin{enumerate}
\item[(i)] $(M,g)$ is of type $(1,1)$ and it is a strict nearly
K\"ahler manifold, \item[(ii)] $(M,g)=(S^6,g_{can})$ and is of
type $(8,8).$
\end{enumerate}
Conversely, if $(M,g)$ is a compact simply-connected strict nearly
K\"ahler $6$-manifold of non-constant curvature then $M$ is of
type $(1,1).$
\end{theorem}

Compact strict nearly K\"ahler manifolds with isometries were
investigated in \cite{MNS05} where it was shown that

\begin{theorem}Let $(M,J,g,\omega_g)$ be a
compact strict nearly K\"ahler $6$-manifold.  If $M$ admits a unit
Killing vector field, then up to finite cover $M$ is isometric to
$S^3\times S^3$ with its standard nearly K\"ahler structure.
\end{theorem}

\begin{remark}\label{holonomy.G2.metric}
The first example of a non-trivial $G_2$ holonomy metric was found
by Bryant \cite{Bry87}, who observed that a cone on the complex
flag manifold $U(3)/T^3$ carries an incomplete metric with
$G_2$-holonomy. The flag ${\rm U}(3)/T^3$ is the twistor space of
the complex projective plane $\bbc\bbp^2$ and as such it also has
a strict nearly K\"ahler structure. As explained in this section,
this therefore is just one possible example. One gets such
non-trivial metrics also for the cones with bases $\bbc\bbp^3$ and
$S^3\times S^3$ with their homogeneous strict nearly K\"ahler
structures. Interestingly, in some cases there exist complete
metrics with $G_2$ holonomy which are smooth deformations of the
asymptotically conical ones. This fact was noticed by Bryant and
Salamon \cite{BrSa89} who constructed complete examples of $G_2$
holonomy metrics on bundles of self-dual $2$-forms over
$\bbc\bbp^2$ and $S^4.$ Replacing the base with any positive QK
orbifold $\calo$ gives complete (in the orbifold sense) metrics on
orbibundles of self-dual $2$-forms over $\calo$. Locally some of
these metrics were considered in \cite{San03}. More complete
examples of explicit $G_2$ holonomy metrics on non-compact
manifolds were obtained by Salamon \cite{Sal04}. $G_2$ holonomy
manifolds with isometric circle actions were investigated by
Apostolov and Salamon \cite{ApSa04}. The first compact examples
are due to the ground breaking work of Joyce \cite{Joy96a}.
\end{remark}

\medskip
\section{Geometries Associated with Sasaki-Einstein
$5$-manifolds}
\medskip

 Like $3$-Sasakian manifolds \Se $5$-manifolds are
naturally associated to other geometries introduced in the
previous section. Of course, each such space $(M^5,\cals)$ comes
with its Calabi-Yau cone $(C(M),\bar{g})$ and, if the \Se
structure $\cals$ is quasi-regular, with its quotient log del
Pezzo surface $(\calz,h).$ But as it turns out, there are {\it
two} more Einstein metrics associated to $g$. The examples of this
section also illustrate how the Theorem \ref{Bar.main.thm1} and
B\"ar's correspondence break down when $(M,g)$ is a manifold
with Killing spinors which is, however, {\it not complete}.

We begin by describing a relation between $5$-dimensional \Se
structures and six-dimensional nearly K\"ahler structures which
was uncovered recently in \cite{FIMU06}. This relation involves
the sine-cones of Definition \ref{conemetric}. We use the notation
$\bar{g}_{s}$ to distinguish the sine-cone metric from the usual
Riemannian cone metric $\bar{g}$. Of course this metric is not
complete, but one can compactify $M$ obtaining a very tractable
stratified space $\bar{M}=N\times [0,\pi]$ with conical
singularities at $t=0$ and $t=\pi.$ Observe the following simple
fact which shows that the Riemannian cone on a sine cone is always
a Riemannian product.

\begin{lemma}\label{CsCm.lemma}\index{sine-cone lemma} Let $(M,g)$ be a Riemannian manifold. Then
the product metric $ds^2=dx^2+dy^2+y^2g$ on $\bbr\times C(M)$ can
be identified with the iterated cone metric on $C(C_s(M))$.
\end{lemma}
\begin{proof}
Consider the map $\bbr^+\times(0,\pi)\ra{1.2}\bbr\times\bbr^+$
given by polar coordinate change $(r,t)\mapsto(x,y)=(r\cos t,r\sin
t),$ where $r>0$ and $t\in(0,\pi).$ We get
$$ds^2=dx^2+dy^2+y^2g=dr^2+r^2dt^2+r^2\sin^2tg=dr^2+r^2(dt^2+\sin^2
tg)\, . \qedhere$$
\end{proof}

So the iterated Riemannian cone $(C(C_s(M)),ds^2)$ has reducible
holonomy $1\times {\rm Hol}(C(M)).$ This leads to

\begin{corollary}\label{sine-coneEinsteinSE}
Let $(N,g)$ be a \Se manifold of dimension $2n+1.$ Then the
sine-cone $C_s(N)$ with the metric $\bar{g}_s=dr^2+(\sin^2 r)g$
is Einstein with Einstein constant $2n+1.$
\end{corollary}

We are particularly interested in the case $n=2.$
Compare Lemma \ref{CsCm.lemma}
with the following result in \cite{Joy00}, Propositions 11.1.1-2:
\begin{proposition}\label{Joy.G2} Let $(M^4,g_4)$ and $(M^6,g_6)$ be Calabi-Yau
manifolds. Let $(\bbr^3, ds^2=dx^2+dy^2+dz^2)$ and $(\bbr,
ds^2=dx^2)$ be the Euclidean spaces. Then
\begin{enumerate}
\item $(\bbr^3\times M^4,g=ds^2+g_4)$ has a natural $G_2$
structure and $g$ has holonomy ${\rm Hol}(g)\subset \BOne_3\times
SU(2)\subset G_2,$ \item $(\bbr\times M^6,g=ds^2+g_6)$ has a
natural $G_2$ structure and $g$ has holonomy ${\rm Hol}(g)\subset
1\times SU(3)\subset G_2.$
\end{enumerate}
\end{proposition}

As long as  $(M^4,g_4)$ and $(M^6,g_6)$ are simply connected then
the products $\bbr^3\times M^4$ and $\bbr\times M^6$ are simply
connected $G_2$-holonomy manifolds with reducible holonomy groups
and parallel Killing spinors. Note that this does not violate
Theorem \ref{Bar.main.thm1} as these spaces are not Riemannian
cones over complete Riemannian manifolds. Using (ii) of Proposition \ref{Joy.G2} we obtain the
following corollary of Theorem \ref{sine-coneEinstein} first obtained in \cite{FIMU06}

\begin{corollary}\label{Se-nearlyK}
Let $(N^5,g)$ be a Sasaki-Einstein manifold. Then the sine cone
$C_s(N^5)=N^5\times (0,\pi)$ with metric $\bar{g}_{s}$ is nearly
K\"ahler of Einstein constant $\lambda=5.$ Furthermore $\bar{g}_s$
approximates pure $SU(3)$ holonomy metric near the cone points.
\end{corollary}

Using Corollary \ref{Se-nearlyK} we obtain a host of examples of
nearly K\"ahler $6$-manifolds with conical singularities by
choosing $N^5$ to be any of the \Se manifolds constructed in
\cite{BGN03c,BGN02b,BG03,Kol04,Kol05b,GMSW04a,GMSW04b,CLPP05,FOW06,CFO07}.
For example, in this way we obtain nearly-K\"ahler metrics on
$N\times (0,\pi)$ where $N$ is any Smale manifold with a \Se
metric such as $S^5$ or $k(S^2\times S^3),$ etc. Note that every
simply connected strict nearly K\"ahler manifold has exactly two
real Killing spinors. So as long as $N^5$ is simply connected
$C_s(N^5)$ will have two real Killing spinors. Using Theorem
\ref{sine-coneEinstein} the \Se metrics constructed in
\cite{BG00a,BGK05,BGKT05,GhKo05,BG06b} in all odd dimensions also
give new Einstein metrics on $C_s(N^{2n+1}).$ For example, one
obtains many positive Einstein metrics on $\grS^{2n+1}\times
(0,\pi)$ where $\grS^{2n+1}$ is any odd dimensional homotopy
sphere bounding a parallelizable manifold. Of course, there are no
Killing spinors unless $n=2.$ Returning to the case of dimension
6, a somewhat more general converse has been obtained in
\cite{FIMU06}, namely

\begin{theorem}\label{nK-Se}
Any totally geodesic hypersurface $N^5$ of a nearly K\"ahler
$6$-manifold $M^6$ admits a \Se structure.
\end{theorem}

The method in \cite{FIMU06} uses the recently developed notion of
hypo $SU(2)$ structure due to Conti and Salamon \cite{CoSa06}. The
study of sine cones appears to have originated in the physics
literature \cite{BiMe03,ADHL03}, but in one dimension higher. Now
recall the following result of Joyce (cf. \cite{Joy00},
Propositions 13.1.2-3)
\begin{proposition}\label{Joy.Spin7}
Let $(M^6,g_6)$ and $(M^7,g_7)$ be Calabi-Yau and $G_2$-holonomy
manifolds, respectively. Let $(\bbr^2, ds^2=dx^2+dy^2)$ and
$(\bbr, ds^2=dx^2)$ be Euclidean spaces. Then
\begin{enumerate}
\item $(\bbr^2\times M^6,g=ds^2+g_6)$ has a natural $Spin(7)$
structure and $g$ has holonomy ${\rm Hol}(g)\subset \BOne_2\times
SU(3)\subset Spin(7),$ \item $(\bbr\times M^7,g=ds^2+g_7)$ has a
natural $Spin(7)$ structure and $g$ has holonomy ${\rm
Hol}(g)\subset 1\times G_2\subset Spin(7).$
\end{enumerate}
\end{proposition}
Again, if $(M^6,g_6)$ and $(M^7,g_7)$ are simply connected so are
the $Spin(7)$-manifolds $\bbr^2\times M^6$ and $\bbr\times M^7$ so
that they have parallel spinors. Not surprisingly, in view of
Lemma \ref{CsCm.lemma} and Proposition \ref{Joy.Spin7}, the sine
cone construction now relates strict nearly K\"ahler geometry in
dimension $6$ to nearly parallel $G_2$ geometry in dimension $7$.
More precisely \cite{BiMe03}

\begin{theorem}\label{NK.to.G2} Let $(N^6,g)$ be a strict nearly K\"ahler
$6$-manifold such that $g$ has Einstein constant $\lambda_6=5$.
Then the manifold $C_s(N)=N^6\times (0,\pi)$ with its sine cone
metric $\bar{g}_{s}$ has a nearly parallel $G_2$ structure with Einstein
constant $\lambda_7=6$ and it approximates pure $G_2$ holonomy
metric near the cone points.
\end{theorem}

\begin{proof} Just as before, starting with $(N^6,g_6)$ we consider
its metric cone $C(N^6)$ with the metric $\bar{g}=dy^2+y^2g_6$ and
the product metric $g_8$ on $\bbr\times C(N^6).$ With the above
choice of the Einstein constant we see that $g_8=dx^2+dy^2+y^2g_6$
must have holonomy ${\rm Hol}(g_8)\subset1\times G_2\subset
Spin(7).$ By Lemma \ref{CsCm.lemma} $g_8$ is a metric cone on the
metric $g_7=dt^2+\sin^2 t g_6,$ which must, therefore, have weak
$G_2$ holonomy and the Einstein constant $\lambda_7=6.$
\end{proof}

Again, any simply connected weak $G_2$-manifold has at least one
Killing spinor. That real Killing spinor on $C_s(N^6)$ will lift
to a parallel spinor on $C(C_s(N^6))=\bbr\times C(N^6)$ which is a
non-complete $Spin(7)$-manifold of holonomy inside $1\times G_2.$
One can iterate the two cases by starting with a compact \Se
$5$-manifold $N^5$ and construct either the cone on the sine cone
of $N^5$ or the sine cone on the sine cone of $N^5$ to obtain a
nearly parallel $G_2$ manifold. We list the Riemannian manifolds
coming from this construction that are irreducible.

\begin{proposition}\label{Se-NK-G2-Spin7}
Let $(N^5,g_5)$ be a  compact \Se manifold which is not of
constant curvature. Then the following have irreducible holonomy
groups:
\begin{enumerate}
\item the manifold $C(N^5)$ with the metric $g_6=dt^2+t^2g_5$ has
holonomy $SU(3);$ \item the manifold $C_s(N^5)=N^5\times (0,\pi)$
with metric $g_6=dt^2 +\sin^2t~g_5$ is strict nearly K\"ahler;
\item the manifold $C_s(C_s(N^5))=N^5\times (0,\pi)\times (0,\pi)$
with the metric $g_7=d\gra^2+\sin^2\gra(dt^2 +\sin^2t~g_5)$ has
a nearly parallel $G_2$ structure.
\end{enumerate}
\end{proposition}

In addition we have the reducible cone metrics:
$C(C_s(N^5))=\bbr\times C(N^5)$ has holonomy in $1\times
SU(3)\subset G_2$ and $C(C_s(C_s(N^5)))=\bbr\times
C(C_s(N^5))=\bbr\times\bbr\times C(N^5)$ has holonomy
$\BOne_2\times SU(3)\subset 1\times G_2\subset Spin(7).$ If $N^5$
is simply connected then  $G_5,g_6$ and $g_7$ admit two Killing
spinors. For a generalization involving conformal factors see
\cite{MoOr07}.

\begin{remark}\label{incomplete.G2} Recall Remark
\ref{complete.G2}. Note that when a nearly parallel $G_2$ metric is
not complete then the type I-III classification is no longer
valid. The group $Spin(7)$ has other subgroups than the ones
listed there and we can consider the following inclusions of
(reducible) holonomies
$${\rm Spin}(7)
\supset G_2\times1\supset SU(3)\times\BOne_2\times \supset
SU(2)\times\BOne_3 \supset\BOne_8 \, .$$ According to the
Friedrich-Kath Theorem \ref{killing.spinor.7D} the middle three
cannot occur as holonomies of Riemannian cones of complete
$7$-manifolds with Killing spinors. But as the discussion of this
section shows, they most certainly can occur as holonomy groups of
Riemannian cones of incomplete nearly parallel $G_2$ metrics. These
metrics can be still separated into three types depending on the
holonomy reduction: say the ones that come from strict nearly
K\"ahler manifolds are generically of type $I_s$ while the ones
that come from \Se $5$-manifolds via the iterated sine cone
construction are of type $II_s$ and of type $III_s$ when $H\subset
SU(3)$ is some proper non-trivial subgroup. On the other hand, it
is not clear what is the relation between the holonomy reduction
and the actual number of Killing spinors one gets in each
case.
\end{remark}

\medskip
\section{Geometric Structures on Manifolds and
Supersymmetry} \medskip

The intricate relationship between supersymmetry and geometric
structures on manifolds was recognized along the way the physics
of supersymmetry slowly evolved from its origins: first globally
supersymmetric field theories (70ties) arose, later came
supergravity theory (80ties), which evolved into
superstring\index{superstrings} theory and conformal field theory
(late 80ties and 90ties), and finally into M-Theory and the
supersymmetric branes of today. At every step the ``first" theory
would quickly led to various generalizations creating many
different new ones: so it is as if after discovering plain vanilla
ice cream one would quickly find oneself in an Italian ice cream
parlor confused and unable to decide which flavor was the right
choice for the hot afternoon. This is a confusion that is possibly
good for one's sense of taste, but many physicists believe that
there should be just one theory, the Grand Unified Theory which
describes our world at any level.\footnote{Actually, string theory
of today appears to offer a rather vast range of vacua (or
possible universes). Such possible predictions have been nicknamed
the {\it string landscape} \cite{Sus03}. This fact has been seen
as a drawback by some, but not all, physicists (see more recent
discussion on {\it landscape} and {\it swampland} in
\cite{Vaf05,OgVa06}). The insistence that the universe we
experience, and this on such a limited scale at best, {\it is the
only Universe}, is largely a matter of `philosophical attitude'
towards science. See the recent book of Leonard Susskind on the
anthropic principle, string theory and the cosmic landscape
\cite{Sus05}.} An interesting way out of this conundrum is to
suggest that even if two theories appear to be completely
different, if both are consistent and admissible, they actually
{\it do} describe the same physical world and, therefore, {\it
they should be dual} to one another in a certain sense. This gave
rise to various duality conjectures such as the Mirror Symmetry
Conjecture or the AdS/CFT Duality Conjecture.

The first observation of how supersymmetry can restrict the
underlying geometry was due to Zumino \cite{Zum79} who discovered
that globally $N=1$ supersymmetric $\sigma$-models in $d=4$
dimensions require that the bosonic fields (particles) of the
theory are local coordinates on a K\"ahler manifold. Later
Alvarez-Gaum\'e and Friedman observed that $N=2$ supersymmetry
requires that the $\sigma$-model manifold be not just K\"ahler but
hyperk\"ahler \cite{AlFr81}. This relation between globally
supersymmetric $\sigma$-models and complex manifolds was used by
Lindstr\"om and Ro\v cek to discover the hyperk\"ahler quotient
construction in \cite{LiRo83, HKLR}.

The late seventies witnessed a series of attempts to incorporate
gravity into the picture which quickly led to the discovery of
various supergravity theories. Again the $N=1$ supergravity-matter
couplings in $d=4$ dimensions require bosonic matter fields to be
coordinates on a K\"ahler manifold with some special properties
\cite{WiBa82} while $N=2$ supergravity demands that the
$\sigma$-model manifold be quaternionic K\"ahler \cite{BaWi83}.
The quaternionic underpinnings of the matter couplings in
supergravity theories lead to the discovery of quaternionic
K\"ahler reduction in \cite{Gal87a, GaLa88}.

At the same time manifolds with Killing spinors emerged as
important players in the physics of the supergravity theory which
in $D=11$ dimensions was first predicted by Nahm \cite{Nah78} and
later constructed by Cremmer, Julia and Scherk \cite{CJS78}. The
well-known Kaluza-Klein trick applied to a $D=11$ supergravity
model is a way of constructing various limiting {\it
compactifications} which would better describe the apparently
four-dimensional physical world we observe. The geometry of such a
compactification is simply a Cartesian product $\bbr^{3,1}\times
M^7$, where $\bbr^{3,1}$ is the Minkowski space-time (or some
other Lorentzian $4$-manifold) and $M^7$ is a compact manifold with
so small a radius that its presence can only be felt and observed
at the quantum level. Many various models for $M^7$ were studied
in the late seventies which by the eighties had already accrued
into a vast physics literature (cf. the extensive three-volume
monograph by Castellani, D'Auria and Fr\'e \cite{CDF91}). Most of
the models assumed a homogeneous space structure on $M^7=G/H$ (see
Chapter V.6 in \cite{CDF91}, for examples). Two things were of key
importance in terms of the required physical properties of the
compactified theory. First, the compact space $M^7$, as a
Riemannian manifold, had to be Einstein of positive scalar
curvature. Second, although one could consider any compact
Einstein space for the compactification, the new theory would no
longer be supersymmetric unless $(M^7,g)$ admitted Killing spinor
fields, and the number of them would be exactly the number of {\it
residual} supersymmetries of the compactified theory. For that
reason compactification models involving $(S^7,g_0)$ were quite
special as they gave the maximally supersymmetric model. However,
early on it was realized that there are other, even homogeneous,
$7$-manifolds of interest. The $Sp(2)$-invariant Jensen metric on
$S^7$, or as physicists correctly nicknamed it, the {\it squashed
$7$-sphere} is one of the examples. Indeed, Jensen's metric admits
exactly one Killing spinor field since it has a nearly parallel $G_2$ structure.
Of course, any of the Einstein geometries in the table of Theorem
\ref{Bar.main.thm1} can be used to obtain such supersymmetric
models.

The $D=11$ supergravity theory only briefly looked liked it was
the Grand Theory of Einstein's dream. It was soon realized that
there are difficulties with getting from $D=11$ supergravity to
the standard model. The theory which was to solve these and other
problems was Superstring Theory and later M-Theory (which is yet
to be constructed). With the arrival of superstring theory and
M-theory, supersymmetry continues its truly remarkable influence
on many different areas of mathematics and physics: from geometry
to analysis and number theory. For instance, once again five, six,
and seven-dimensional manifolds admitting real Killing spinors
have become of interest because of the so called AdS/CFT Duality.
Such manifolds have emerged naturally in the context of $p$-brane
solutions in superstring theory. These so-called $p$-branes,
``near the horizon" are modelled by the pseudo-Riemannian geometry
of the product ${\rm AdS}_{p+2}\times M$, where ${\rm AdS}_{p+2}$
is the $(p+2)$-dimensional anti-de-Sitter space (a Lorentzian
version of a space of constant sectional curvature) and $(M, g)$
is a Riemannian manifold of dimension $d=D-p-2$. Here $D$ is the
dimension of the original supersymmetric theory. In the most
interesting cases of M2-branes, M5-branes, and D3-branes $D$
equals either 11 (M$p$-branes of M-theory) or 10 (D$p$-branes in
type IIA or type IIB string theory). String theorists are
particularly interested in those vacua of the form ${\rm
AdS}_{p+2}\times M$ that preserve some residual supersymmetry. It
turns out that this requirement imposes constraints on the
geometry of the Einstein manifold $M$ which is forced to admit
real Killing spinors. Depending on the dimension $d$, the possible
geometries of $M$ are as follows:

\begin{table}[h!]
\centering
\begin{tabular}{|c|c|c|}\hline
d & Geometry of M & $(\mu, \Bar\mu)$\\
\hline \hline
any & round sphere & $(1,1)$\\
$7$ & nearly parallel $G_2$ & $(\frac18,0)$\\
  & Sasaki--Einstein & $(\frac14,0)$\\
  & 3-Sasakian & $(\frac38, 0)$\\
$6$ & nearly K\"ahler & $(\frac18,\frac18)$\\
$5$ & Sasaki--Einstein & $(\frac14,\tfrac14)$\\
\hline
\end{tabular}
\end{table}
\noindent where the notation $(\mu,\bar{\mu}),$ which is common in
the physics literature, represents the ratio of the number of real
Killing spinors of type $(p,q)$ to the maximal number of real
Killing spinors that can occur in the given dimension. This
maximum is, of course, realized by the round sphere of that
dimension. So this table is just a translation of the table of
Theorem \ref{Bar.main.thm1} for the special dimensions that occur
in the models used by the physicists.

Furthermore, given a $p$-brane solution of the above type, the
interpolation between ${\rm AdS}_{p+2}\times M$ and
$\bbr^{p,1}\times \calc(M)$ leads to a conjectured duality between
the supersymmetric background of the form ${\rm AdS}_{p+2}\times M$
and a $(p+1)$-dimensional superconformal field theory of $n$
coincident $p$-branes located at the conical singularity of the
$\bbr^{p,1}\times \calc(M)$ vacuum. This is a generalized version of
the Maldacena or AdS/CFT Conjecture \cite{Mal99}\index{AdS/CFT
Duality}. In the case of D3-branes of string theory the relevant
near horizon geometry is that of ${\rm AdS}_{5}\times M$, where $M$
is a \Se 5-manifold. The D3-brane solution interpolates between
${\rm AdS}_{5}\times M$ and $\bbr^{3,1}\times \calc(M)$, where the
cone $\calc(M)$ is a Calabi-Yau threefold. In its original version
the Maldacena conjecture (also known as AdS/CFT duality) states that
the 't Hooft large $n$ limit of $N=4$ supersymmetric Yang-Mills
theory with gauge group $SU(n)$ is dual to type IIB superstring
theory on ${\rm AdS}_{5}\times S^5$ \cite{Mal99}. This conjecture
was further examined by Klebanov and Witten \cite{KlWi99} for the
type IIB theory on ${\rm AdS}_{5}\times T^{1,1}$, where $T^{1,1}$ is
the other homogeneous \Se $5$-manifold $T^{1,1}=S^2\times S^3$ and
the Calabi-Yau 3-fold $\calc(T^{1,1})$ is simply the quadric cone in
$\bbc^4$. Using the well-known fact that $\calc(T^{1,1})$ is a
K\"ahler quotient of $\bbc^4$ (or, equivalently, that $S^2\times
S^3$ is a Sasaki-Einstein quotient of $S^7$), a dual super
Yang-Mills theory was proposed, representing D3-branes at the
conical singularities. In the framework of D3-branes and the AdS/CFT
duality the question of what are all the possible near horizon
geometries $M$ and $\calc(M)$ might be of importance. Much of the
interest in \Se manifolds is precisely due to the fact that each
such explicit metric, among other things, provides a useful model to
test the AdS/CFT duality. We refer the reader interested in the
mathematics and physics of the AdS/CFT duality to the recent book in
the same series \cite{Biq05}. In particular, in this context,
Sasaki-Einstein geometry  is discussed in one of the articles there
\cite{GMSW05pc}.

\begin{remark}{\bf [$G_2$ holonomy manifolds unification scale and proton decay]}
Until quite recently the interest in $7$-manifolds with $G_2$
holonomy as a source of possible physical models was tempered by
the fact the Kaluza-Klein compactifications on smooth and complete
manifolds of this type led to models with no charged particles.
All this has dramatically changed in the last few years largely
because of some new developments in M-theory. Perhaps the most
compelling reasons for reconsidering such $7$-manifolds was
offered by Atiyah and Witten who considered the dynamics on
manifolds with $G_2$ holonomy which are asymptotically conical
\cite{AtWi02}. The three models of cones on the homogeneous nearly
K\"ahler manifolds mentioned earlier are of particular interest,
but Atiyah and Witten consider other cases which include orbifold
(quotient) singularities. Among other things they point to a very
interesting connection between Kronheimer's quotient construction
of the ALE metrics \cite{Kro89a, Kro89b} and asymptotically
conical manifolds with $G_2$-holonomy.  To explain the connection,
consider Kronheimer's construction for $\Gamma=\bbz_{n+1}.$
Suppose one chooses a circle $S^1_{k,l}\simeq{\rm U}(1)\in
K(\bbz_{n+1})={\rm U}(1)^n$ and then one considers a $7$-manifold
obtained by performing Kronheimer's HK quotient construction with
zero momentum level $(\boldsymbol{\xi}=0)$ while "forgetting" the
three moment map equations corresponding to this particular
circle. An equivalent way of looking at this situation is to take
the Kronheimer quotient with nonzero momentum
$\boldsymbol{\xi}=a\in \gs\gp(1)$ but only for the moment map of
the chosen circle $S^1_{k,l}$ (such $\boldsymbol{\xi}$ is never in
the "good set") and then consider the fibration of singular
Kronheimer quotients over a $3$-dimensional base parameter space.
Algebraically this corresponds to a partial resolution of the
quotient singularity and this resolution depends on the choice of
$S^1_{k,l}$, hence $\boldsymbol{\xi}.$ This example was first
introduced in \cite{AtWi02}. It can be shown that the $7$-manifold
is actually a cone on the complex weighted projective $3$-space
with weights $(k,k,l,l),$ where $k+l=n+1.$ It then follows from
the physical model considered that such a cone should admit a
metric with $G_2$ holonomy. However, unlike the homogeneous cones
over the four homogeneous strict nearly K\"ahler manifolds of
Theorem \ref{But.thm}, the metric in this case is not known
explicitly. This construction appears to differ from all previous
geometric constructions of metrics with $G_2$ holonomy. One can
consider similar constructions for other choices of $S^1\subset
K(\Gamma)$ \cite{BeBr02}.

In \cite{FrWi03} using a specific models of M-theory
compactifications on manifolds with $G_2$ holonomy, Friedman and
Witten address the fundamental questions concerning the
unification scale ({\it i.e.}, the scale at which the Standard
Model of $SU(3)\times SU(2)\times U(1)$ unifies in a single gauge
group) and proton decay. The authors point out that the results
obtained are model dependent, but some of the calculations and
conclusions apply to a variety of different models.
\end{remark}

\newcommand{\etalchar}[1]{$^{#1}$}
\def\cprime{$'$} \def\cprime{$'$} \def\cprime{$'$} \def\cprime{$'$}
  \def\cprime{$'$} \def\cprime{$'$} \def\cprime{$'$} \def\cprime{$'$}
  \def\cdprime{$''$}
\providecommand{\bysame}{\leavevmode\hbox to3em{\hrulefill}\thinspace}
\providecommand{\MR}{\relax\ifhmode\unskip\space\fi MR }
% \MRhref is called by the amsart/book/proc definition of \MR.
\providecommand{\MRhref}[2]{%
  \href{http://www.ams.org/mathscinet-getitem?mr=#1}{#2}
}
\providecommand{\href}[2]{#2}

\printindex

\end{document}